\documentclass[a4paper,10pt]{article}
\usepackage[utf8]{inputenc}
\usepackage{amsmath,amsthm}

\usepackage{stmaryrd} 
\usepackage{dsfont} 
\usepackage{graphicx}
\usepackage{enumerate}
\usepackage{subfigure}

\usepackage{xcolor}

\theoremstyle{definition}
\newtheorem{definition}{Definition}[]
\theoremstyle{definition}
\newtheorem{remark}{Remark}[]

\newcommand{\ddt}[1]{\frac{d #1}{dt}}

\newcommand{\inproduct}[2]{\left\langle #1 , #2 \right\rangle}
\newcommand{\jump}[1]{\left\llbracket #1 \right\rrbracket}

\newcommand{\divergence}[1]{\nabla \cdot #1}

\newcommand{\transpose}{{\sf T}}
\newcommand{\norm}[1]{\lVert #1 \rVert} 
\DeclareMathOperator{\gr}{gr} 

\newcommand{\RR}{\mathds{R}} 

\DeclareMathOperator{\cfl}{CFL} 
\DeclareMathOperator{\thetaref}{\theta_{ref}}
\DeclareMathOperator{\thetanum}{\theta_{num}}
\DeclareMathOperator{\xclnum}{x_{cl,num}}
\DeclareMathOperator{\xclref}{x_{cl,ref}}
\DeclareMathOperator{\errortheta}{\mathcal{E}_\theta^\infty}
\DeclareMathOperator{\errorcl}{\mathcal{E}_{cl}^\infty}

\newcommand{\move}{\mathcal{M}} 

\newcommand{\nsigma}{n_\Sigma}

\newcommand{\ndomega}{n_{\partial\Omega}}
\newcommand{\ngamma}{n_\Gamma}
\newcommand{\tgamma}{t_\Gamma}

\newcommand{\pdomega}{\mathcal{P}_{\partial\Omega}}

\newcommand{\normalspeed}{V_\Sigma} 

\usepackage{authblk} 
\usepackage{hyperref}
\usepackage{doi}
\usepackage{cite}
\usepackage{microtype}
\usepackage[a4paper,left=2cm,right=2cm,top=3cm,bottom=3cm]{geometry}

\title{Contact line advection using the geometrical\\ Volume-of-Fluid method}
\author[1]{Mathis Fricke\thanks{fricke@mma.tu-darmstadt.de}}
\author[1]{Tomislav Mari\'{c}\thanks{maric@mma.tu-darmstadt.de}}
\author[1]{Dieter Bothe\thanks{bothe@mma.tu-darmstadt.de}}

\affil[1]{Mathematical Modeling and Analysis Group, TU Darmstadt, Alarich-Weiss-Straße 10,\newline 64287 Darmstadt, Germany}

\begin{document}

\maketitle

\begin{abstract}
We consider the interface advection problem  by a prescribed velocity field in the special case when the interface intersects the domain boundary, i.e.\ in the presence of a contact line. This problem emerges from the discretization of continuum models for dynamic wetting. The \emph{kinematic evolution equation} for the dynamic contact angle (Fricke et al., 2019) expresses the fundamental relationship between the rate of change of the contact angle and the structure of the transporting velocity field. The goal of the present work is to develop an interface advection method that is consistent with the fundamental kinematics and transports the contact angle correctly with respect to a prescribed velocity field. In order to verify the advection method, the kinematic evolution equation is solved numerically and analytically (for special cases).\\
We employ the geometrical Volume-of-Fluid (VOF) method on a structured Cartesian grid to solve the hyperbolic transport equation for the interface in two spatial dimensions. We introduce generalizations of the Youngs and ELVIRA methods to reconstruct the interface close to the domain boundary. Both methods deliver first-order convergent results for the motion of the contact line. However, the Boundary Youngs method shows strong oscillations in the numerical contact angle that do not converge with mesh refinement. In contrast to that, the Boundary ELVIRA method provides linear convergence of the numerical contact angle transport.
\end{abstract}

This preprint was submitted and accepted for publication in the \emph{Journal of Computational Physics}.\newline When citing this work, please refer to the journal article: \textbf{DOI:} \href{https://doi.org/10.1016/j.jcp.2019.109221}{10.1016/j.jcp.2019.109221}.\newline
\newline
\textbf{Keywords:} Volume-of-Fluid, Interface Reconstruction, Dynamic Contact Angle, Kinematics



\section{Introduction}
\begin{figure}[ht]
\begin{minipage}{0.45\columnwidth}
\includegraphics[height=3cm]{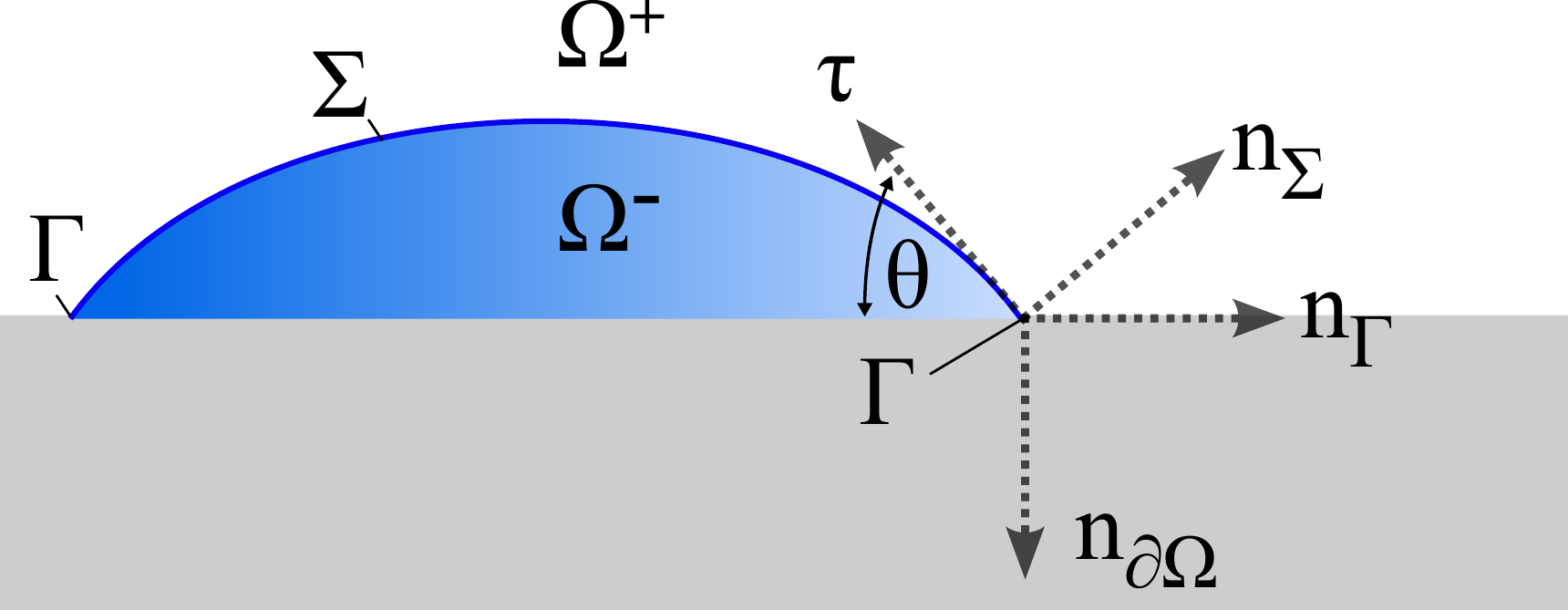}
\caption{Notation.}
\label{fig:notation}
\end{minipage}
\hfil
\begin{minipage}{0.45\columnwidth}
\includegraphics[height=3cm]{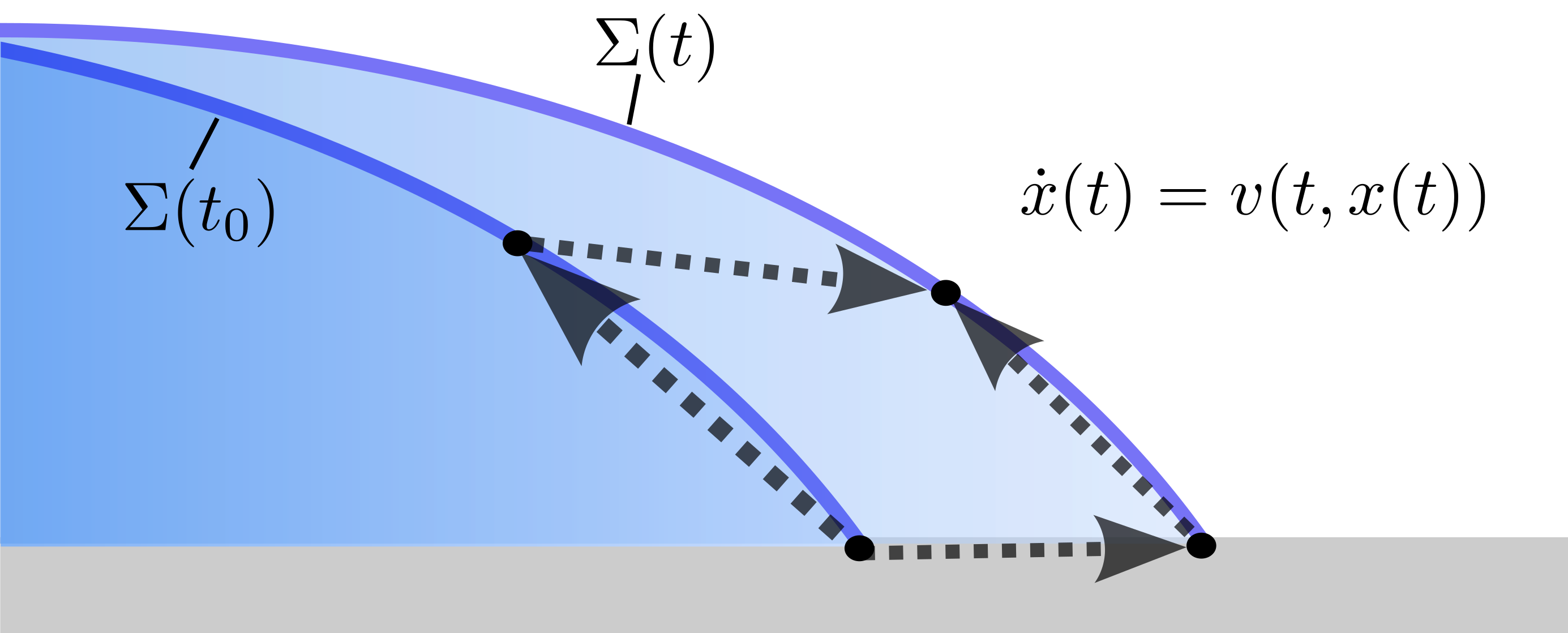}
\caption{Kinematic transport of the contact angle.}
\label{fig:kinematic_transport}
\end{minipage}
\end{figure}

The present article deals with the passive transport of an interface by a prescribed velocity field, a well-known problem in the numerical description of multiphase flows called \emph{interface advection} (see, e.g., \cite{Rider1998,Tryggvason2011}). In particular, we consider the special case when the interface intersects the domain boundary and a so-called \emph{contact line} is formed. This problem arises, for example, as part of the mathematical description of \emph{wetting} processes, where a liquid displaces another liquid (or a gas) in order to wet the surface of a third solid phase. In the latter case, the velocity field is a solution of the two-phase Navier Stokes equations with appropriate boundary and transmission conditions modeling the physics of wetting.\\
\\
In the present work, we focus on the \emph{kinematics} of wetting, hence treating the velocity field as \emph{given}, while keeping in mind that it is a solution to some continuum mechanical model. The continuum mechanical modeling of wetting is a separate issue that is not addressed here (see, e.g., \cite{Bonn.2009,Shikhmurzaev.2008}). The kinematics of wetting is studied analytically in \cite{Fricke.2019, Fricke.2018} where it is shown that the contact line advection problem is a well-posed initial value problem if the velocity field is sufficiently regular and tangential to the domain boundary. Hence the full dynamics of the interface can be inferred from the (time-dependent) velocity field and the initial interface configuration. Note that this is true even though in the full continuum mechanical model the interface shape is typically strongly coupled to the flow. The reason is that the velocity field contains enough information to reconstruct the evolution of the interface from the initial configuration (see Figure~\ref{fig:kinematic_transport} for a sketch of the idea). In particular, the motion of the contact line as well as the evolution of the contact angle, i.e.\ the angle of intersection between the fluid-fluid interface and the solid boundary (see Figure~\ref{fig:notation}), can be computed from the knowledge of the velocity field and the initial geometry. In fact, the rate of change of the contact angle $\dot\theta$ satisfies the kinematic evolution equation
\begin{align}
\label{eqn:theta_evolution_intro}
\dot\theta = \frac{\partial v}{\partial \tau} \cdot \nsigma 
\end{align}
derived in \cite{Fricke.2019} (see Figure~\ref{fig:notation}). The goal of the present work is to develop an advection scheme for the geometrical VOF method that allows to transport the contact angle consistently with the fundamental kinematic relation \eqref{eqn:theta_evolution_intro}, which will serve as a reference to validate the method (see Section~\ref{section:kinematics}).\\
\\
On the other hand, mathematical models for dynamic wetting usually \emph{prescribe} a boundary condition for the contact angle. The contact angle is an important physical parameter characterizing the wettability of the solid surface. Mathematically, it is typically prescribed as a fixed value or as a function of the speed of the contact line (see, e.g., \cite{Gennes.1985}). Note that in view of \eqref{eqn:theta_evolution_intro}, a boundary condition for the contact angle leads to a consistency condition for the velocity field. In order to avoid the moving contact line paradox \cite{Huh.1971}, tangential slip is usually introduced on the solid boundary, at least close to the contact line.\\
\\
The contact angle boundary condition also plays an important role in the numerical simulation of moving contact lines (see \cite{Sui.2014} for a recent review). In the present study, we employ the Volume-of-Fluid (VOF) method (see \cite{Hirt1981, Rider1998}) to numerically track the location of the interface. The VOF method has been successfully adapted for the simulation of moving contact lines by several authors (see, e.g., \cite{Afkhami2009,Afkhami.2009b,Afkhami2018,Dupont.2010,Fath2015,Kunkelmann2009,Raeini2012,Renardy.2001,Sikalo2005}). Since the advection of the volume fraction field is mostly discretized explicitly in time, the contact angle boundary condition is usually not satisfied after a transport step. Instead, the contact angle is enforced by an explicit adjustment of the interface orientation at computational cells located at the boundary. Enforcing the contact angle as described above is, however, not consistent with the kinematics of moving contact lines as discussed in \cite{Fricke.2019},\cite{Fricke.2018}. We observed that it can also be a source for numerical instabilities in the vicinity of the contact line (see \cite{Raeini2012} for a discussion of instabilities close to the contact line).\\

\paragraph{Notation:} We consider a moving $\mathcal{C}^{1,2}$-hypersurface $\{\overline\Sigma(t)\}_{t \in I}$ with boundary (as defined in the appendix) embedded in two or three-dimensional Euclidean space. Each instantaneous interface $\overline\Sigma(t)$ is assumed to be an orientable $\mathcal{C}^2$-hypersurface with unit normal vector field $\nsigma(t,\cdot)$ and boundary $\partial\Sigma(t)$ that is contained in the planar domain boundary denoted by $\partial\Omega$. The contact line at time $t$ is denoted by
\[ \Gamma(t) = \partial\Sigma(t) \subset \partial\Omega.  \]
Given a point $x \in \Gamma(t)$, the contact angle $\theta \in (0,\pi)$ is defined by the relation
\begin{align}
\label{eqn:theta_definition}
\cos \theta(t,x) = - \inproduct{\nsigma(t,x)}{\ndomega(x)}, 
\end{align}
where $\inproduct{\cdot}{\cdot}$ is the Euclidean scalar product and $\ndomega$ is the outer unit normal field to $\partial\Omega$. In order to describe the motion of the interface, the notion of the \emph{speed of normal displacement} is required. If $\gamma$ is a continuously differentiable curve such that
\[ \gamma(t) \in \overline\Sigma(t) \quad \text{and} \quad \gamma(t_0) = x_0 \in \overline\Sigma(t_0), \]
then the speed of normal displacement at the point $x_0$ at time $t_0$ is given as\footnote{It can be shown that this definition is independent of the choice of the curve $\gamma$, see e.g. \cite{Pruss.2016}.}
\begin{align}
\normalspeed(t_0,x_0) = \inproduct{\gamma'(t_0)}{\nsigma(t_0,x_0)}.
\end{align}
It is convenient to introduce a contact line normal vector (see Figure~\ref{fig:notation}) via projection as
\[ \ngamma = \frac{\pdomega \nsigma}{\norm{\pdomega \nsigma}}, \]
where $\pdomega = \mathds{1} - \ndomega \otimes \ndomega$ is the local orthogonal projection onto $\partial\Omega$. The above expression is well-defined since we only consider the \emph{partial wetting case} characterized by $0 < \theta < \pi$. To construct a local orthonormal basis in three dimensions, one further defines
\[ \tgamma = \ngamma \times \ndomega. \]

\paragraph{Problem formulation:} The motion of a material interface is governed by the kinematic condition
\begin{align}
\label{eqn:kinematic_boundary_condition}
 \normalspeed = \inproduct{v}{\nsigma} \quad \text{on} \quad \gr\overline\Sigma,
\end{align}
where $v$ is the transporting velocity field. Here it is assumed that $v \in \mathcal{C}^1(\overline\Omega)$ is divergence free and tangential to the domain boundary, i.e.
\begin{align}
\divergence{v} = 0 \quad &\text{in} \ \Omega,\label{eqn:zero_divergence}\\
v \cdot \ndomega = 0 \quad &\text{on} \ \partial\Omega.\label{eqn:impermeability}
\end{align}
For simplicity, we further assume that the solid boundary is \emph{planar}. This assumption is not essential and may be dropped to study the contact angle evolution over curved surfaces. Note also that incompressibility of the flow is not necessary to study the contact angle evolution. However, it is assumed here since the VOF method is most commonly applied to incompressible flows.\\
\\
We emphasize that the assumption of a globally continuously differentiable velocity field is \emph{not} met in a typical multiphase flow, where the velocity gradient admits a jump which is controlled by the interfacial transmission condition for the stress. However, it has been shown in \cite{Fricke.2019} Lemma~8, that in the case of two spatial dimensions\footnote{A similar statement holds in three spatial dimensions. For that case, one can show the continuity property $\inproduct{\jump{\nabla v} \alpha}{\beta} = 0$ at $\Gamma$, where $\jump{\cdot}$ denotes the jump over the interface and $\alpha$, $\beta$ are arbitrary vectors in the plane spanned by $\ndomega$ and $\ngamma$, see \cite{Fricke.2019}.} the conditions 
\[ \jump{v} = 0 \ \ \text{on} \ \gr\overline\Sigma, \quad \nabla \cdot v = 0 \ \ \text{in} \ \ \Omega\setminus\Sigma(t), \quad v \cdot \ndomega = 0 \ \ \text{on} \ \ \partial\Omega \] 
imply continuity of $\nabla v$ \emph{at} the contact line (under certain regularity assumptions). For simplicity, we therefore assume a globally continuously differentiable field.\\
\\
In this work, we focus on the Volume-of-Fluid method \cite{Hirt1981,Rider1998} for advecting the interface, whose reconstruction algorithm must be adapted in order to achieve second-order convergence in the near-wall region. However, the contact line advection problem can be used to verify near-wall advection of any other interface advection method. For example, a brief discussion of the contact line advection problem with the Level Set Method \cite{Sethian1996,Sussman1998} in two dimensions can be found in \cite{Fricke.2019b}. An open research data record containing the full C++-implementation of the Level Set Method used in \cite{Fricke.2019b} together with a number of computational examples is available online, see \cite{OnlineRepoPAMM}.

\section[Kinematic evolution equation]{Kinematic evolution equation}
\label{section:kinematics}
The kinematic evolution equation (see \cite{Fricke.2019},\cite{Fricke.2018}) allows to compute the rate of change of the contact angle in terms of the transporting velocity field. The evolution of the contact angle is considered along the trajectories (or ``characteristics'') of the flow, i.e.\ along solutions of the ordinary differential equation
\begin{align}
\label{eqn:characteristics_ode}
\dot{x}(t) = v(t,x(t)), \quad x(t_0) = x_0 \in \overline\Sigma(t_0).
\end{align}
The conditions \eqref{eqn:kinematic_boundary_condition}, \eqref{eqn:zero_divergence} and \eqref{eqn:impermeability} imply that both $\gr\Gamma$ and $\gr\Sigma=\gr\overline\Sigma\setminus\gr\Gamma$ are \emph{invariant} for the flow generated by \eqref{eqn:characteristics_ode}, i.e. a solution of \eqref{eqn:characteristics_ode} starting on the contact line (the fluid-fluid interface) will stay on the contact line (the fluid-fluid interface); see \cite{Fricke.2019} for a proof. Therefore, it is possible to study the time evolution of the contact angle along a solution of \eqref{eqn:characteristics_ode} starting at the contact line. Along such a curve, one can define the \emph{Lagrangian time-derivative} of a quantity $\psi \in \mathcal{C}^1(\gr\overline\Sigma)$ according to
\begin{align*}
\dot{\psi}(t_0,x_0) = \ddt{} \, \psi(t,x(t))_{|t=t_0}.
\end{align*}

Given a $\mathcal{C}^1$-velocity field, the Lagrangian time-derivative of the interface normal vector is given by (see \cite{Fricke.2019}, \cite{Fricke.2018})
\begin{equation}
\label{eqn:ode_nsigma_evolution}
\begin{aligned}
\dot{n}_\Sigma(t,x) &= -[\mathds{1}-\nsigma(t,x)\otimes\nsigma(t,x)] \, \nabla v(t,x)^\transpose \nsigma(t,x).
\end{aligned}
\end{equation}
Clearly, the evolution of the contact angle can be inferred from the normal vector evolution via the relation \eqref{eqn:theta_definition}. Moreover, it has been shown \cite{Fricke.2019}, \cite{Fricke.2018} that, for a planar boundary\footnote{Note that \eqref{eqn:ode_nsigma_evolution} also holds if $\partial\Omega$ is not planar.}, the contact angle follows the evolution equation
\begin{align}
\label{eqn:theta_evolution_equation}
\dot{\theta}(t,x) = \inproduct{\partial_\tau v(t,x)}{\nsigma(t,x)},
\end{align}
where the tangential direction $\tau$ is defined as (see Figure~\ref{fig:notation})
\[ \tau(t,x) = - \cos[\theta(t,x)] \, \ngamma(t,x) - \sin [\theta(t,x)] \, \ndomega(t,x). \]
The goal of the present work is to verify the numerical solution delivered by the VOF method against an analytical (if available) or numerical solution of \eqref{eqn:theta_evolution_equation}.

\subsection[An analytical solution for linear velocity fields in 2D]{An analytical solution for linear velocity fields in 2D}
We first consider the case of general linear divergence free velocity fields in 2D. In this case, the velocity gradient $\nabla v$ is constant in space and time and the ODE system \eqref{eqn:characteristics_ode} and \eqref{eqn:ode_nsigma_evolution} is \emph{explicitly} solvable. Note that this also provides a local approximation for general differentiable velocity fields.\\
\\
We choose a Cartesian coordinate system $(x_1,x_2)$ such that the solid wall is represented by $x_2 = 0$. We consider a velocity field of the form
\begin{align}
v(x_1,x_2) = (v_0 + c_1 x_1 + c_2 x_2, -c_1 x_2).
\end{align}
The coefficients $c_1$ and $c_2$ in this formulation have the dimension of $s^{-1}$. Therefore, it is more convenient to choose  a length scale $L$ and a time scale $T$ and write
\begin{align}
\label{eqn:general_linear_field}
\frac{v(x_1,x_2)}{L/T} = (\hat{v}_0 + \hat{c}_1 \hat{x}_1 + \hat{c}_2 \hat{x}_2, -\hat{c}_1 \hat{x}_2) =: \hat{v}(\hat{x}_1,\hat{x}_2)
\end{align}
with the non-dimensional quantities $\hat{x}_i = x/L$, $\hat{c}_i = c_i T$ and $\hat{v}_0 = (T v_0)/L$. In the following, we will use the formulation \eqref{eqn:general_linear_field} while dropping the hats. For a field of this form, the (constant) gradient is given by
\[ \nabla v = \begin{pmatrix} c_1 & c_2 \\ 0 & -c_1 \end{pmatrix}. \]
\paragraph{Motion of the contact line:} The motion of the contact line is determined by the ordinary differential equation
\[ \dot{x}_1(t) = v_1(x_1(t),0) = v_0 + c_1 x_1(t), \quad x_1(0) = x_1^0. \]
The unique solution of the above initial value problem is
\begin{align}
\label{eqn:contactline_motion_navier}
x_1(t) = x_1^0 e^{c_1 t} + \frac{v_0}{c_1} \left(e^{c_1 t} -1 \right) \quad \text{for} \ c_1 \neq 0
\end{align}
and $x_1(t) = x_1^0 + v_0 t$ for $c_1 = 0$.

\paragraph{Contact angle evolution:} Note that the constancy of $\nabla v$ decouples the system \eqref{eqn:characteristics_ode} and \eqref{eqn:ode_nsigma_evolution}. Hence, the evolution of the normal vector can be solved independently of the evolution of the contact point. To find the solution, we make use of the fact that in two dimensions the normal vector $\nsigma$ is, up to a reflection, uniquely determined by the contact angle $\theta$. Given a contact angle $\theta$, the two possibilities are
\begin{align*}
\nsigma^l = \begin{pmatrix} - \sin \theta \\ \cos \theta \end{pmatrix} \quad \text{and} \quad \nsigma^r = \begin{pmatrix} \sin \theta \\ \cos \theta \end{pmatrix}.
\end{align*}
In the case of a droplet (and for $\theta < \pi$), this corresponds to the two distinct contact points (left and right). The corresponding expressions for $\tau$ are
\begin{align*}
\tau^l = \begin{pmatrix} \cos \theta \\ \sin \theta \end{pmatrix} \quad \text{and} \quad \tau^r = \begin{pmatrix} -\cos \theta \\ \sin \theta \end{pmatrix}.
\end{align*}
This allows to infer the evolution of $\theta$ for the left and the right contact point directly from \eqref{eqn:theta_evolution_equation} without the need to solve the system \eqref{eqn:ode_nsigma_evolution}.\\
\\
Inserting the expressions for $\nabla v$, $\nsigma^{l,r}$ and $\tau^{l,r}$ to \eqref{eqn:theta_evolution_equation} yields the nonlinear ordinary differential equation\footnote{Note that we now use the Lagrangian formulation and write $\theta(t)$ for the contact angle at $x(t) \in \Gamma(t)$ where $x(t)$ is a trajectory of the flow, i.e.\ a solution of \eqref{eqn:characteristics_ode}.}
\begin{align}
\label{eqn:theta_evolution_linear_2d}
\dot{\theta}(t) = \pm c_2 \sin^2 \theta - 2 c_1 \sin \theta \cos \theta, \quad \theta(0) = \theta_0
\end{align}
with the ``$+$'' for the evolution of the right contact point and the ``$-$'' for the evolution of the left contact point. The unique solution of the above initial value problem is given by the formula (see Appendix for a derivation)
\begin{align}
\label{eqn:exact_solution_2d}
\theta(t) = \frac{\pi}{2} + \arctan\left(-\cot \theta_0 \, e^{2c_1 t} \pm c_2 \frac{e^{2c_1 t} -1}{2c_1} \right).
\end{align}

\paragraph{Remark:} Obviously, the solution is independent of the parameter $v_0$. This is to be expected since the two differential equations \eqref{eqn:characteristics_ode} and \eqref{eqn:ode_nsigma_evolution} decouple and the parameter $v_0$ can be eliminated by a change of the frame of reference. Moreover, the evolution of the left and the right contact point is identical if $\theta_0^r = \theta_0^l$ and $c_2 = 0$. Finally, we note that \eqref{eqn:exact_solution_2d} has a well-defined limit for $c_1 \rightarrow 0$ since
\[ \lim_{c_1 \rightarrow 0} \left(-\cot \theta_0 \, e^{2c_1 t} \pm c_2 \frac{e^{2c_1 t} -1}{2c_1}\right) = -\cot \theta_0 \pm c_2 t.  \]

\section{Numerical Method}
The numerical studies in the present work are carried out with the Volume of Fluid in-house code \textit{Free Surface 3D} (FS3D) originally developed by Martin Rieber \cite{Rieber2004}. Since then, FS3D has been further developed at the University of Stuttgart (see, e.g., \cite{Eisenschmidt2016} and more references given there) and at the Technical University of Darmstadt (see, e.g., \cite{Fath2015} and references therein).

\subsection{Geometrical Volume-of-Fluid Method}
The Volume of Fluid method makes use of the phase indicator function for one of the phases $\Omega^\pm(t)$, which are separated by the interface $\Sigma(t)$, i.e.
\begin{align}
     \chi(t,x) = \begin{cases} 1 & \text{if} \quad x \in \Omega^-(t) \\ 0 & \text{if} \quad x \notin \Omega^-(t) \end{cases},
\end{align}
to describe the interface. Formally, the kinematic condition \eqref{eqn:kinematic_boundary_condition} translates to the transport equation
\[  \partial_t \chi + v \cdot \nabla \chi = 0\]
for the indicator function, which is its conservative forms reads as
\begin{align}
\label{eqn:chi_transport_conservative}
\partial_t \chi + \nabla \cdot (v\chi) = \chi (\nabla \cdot v).
\end{align}
Note that the right-hand side vanishes for an incompressible flow. However, it is kept in the discretization to enhance volume conservation properties \cite{Rider1998}. The latter approach is called ``divergence correction''. Integrating \eqref{eqn:chi_transport_conservative} over a control volume in two spatial dimensions leads to
\begin{align}
\label{eqn:continuous_alpha_evolution}
\ddt{} \alpha_{ij}(t) = - \frac{1}{|V_{ij}|} \int_{\partial V_{ij}} \chi \,  v \cdot n \, dA + \frac{1}{|V_{ij}|} \int_{V_{ij}} \chi (\nabla \cdot v) \, dV, 
\end{align}
where the volume fraction $\alpha_{ij}$ associated with the control volume $V_{ij}$ at time $t$ is introduced as
\begin{align}
\alpha_{ij}(t) = \frac{1}{|V_{ij}|} \int_{V_{ij}} \chi \, dV.
\end{align}
A temporal integration of \eqref{eqn:continuous_alpha_evolution} leads to the exact transport equation of the volume fraction field given as
\begin{align}
\label{eqn:continuous_alpha_evolution_v2}
\alpha_{ij}(t+\Delta t) = \alpha_{ij}(t) - \frac{1}{|V_{ij}|} \int_{t}^{t+\Delta t} \int_{\partial V_{ij}} \chi \,  v \cdot n \, dA \, d\tau + \frac{1}{|V_{ij}|} \int_{t}^{t+\Delta t} \int_{V_{ij}} \chi (\nabla \cdot v) \, dV \, d\tau.
\end{align}
The conservative form of \eqref{eqn:continuous_alpha_evolution_v2} allows for the exact conservation of the phase volume also in the discrete case.  Algebraic Volume-of-Fluid methods approximate the solution of \eqref{eqn:continuous_alpha_evolution_v2} using a Finite Volume discretization that relies on an interpolation of $\chi$ at $\partial V_{ij}$. The sharp jump in $\chi$ between phases "$+$" and "$-$" then causes oscillations in the numerical solution, that are counteracted by adding artificial diffusion which, in turn, leads to artificial smearing of the interface \cite{Desphande2012}. \emph{Geometric} Volume-of-Fluid methods are more accurate than their algebraic counterparts, as they approximate the numerical flux, which can be decomposed in a sum over faces $A_f$ of the control volume according to
\[ \int_{t}^{t+\Delta t} \int_{\partial V_{ij}} \chi \,  v \cdot n \, dA \, d\tau = \sum_f \int_{t}^{t+\Delta t} \int_{A_f} \chi \,  v \cdot n \, dA \, d\tau, \]
by reconstructing a sharp geometrical approximation of the indicator function $\chi(t, \cdot)$ and subsequently approximating the integral using geometrical methods.\\
\\
Within the FS3D code, equation \eqref{eqn:continuous_alpha_evolution_v2} is solved on a structured Eulerian grid using an operator splitting method as described below (see also \cite{Rider1998,Rieber2004,Eisenschmidt2016}).

\begin{figure}[ht]
\centering
\includegraphics[width=0.2\columnwidth]{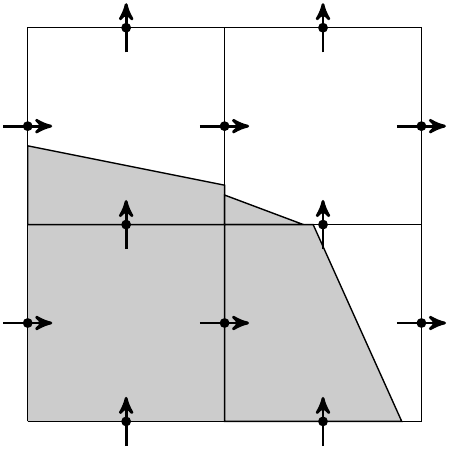}
\caption{Staggered grid for velocity components.}
\label{fig:staggered_grid}
\end{figure}

\paragraph{Implementation in FS3D:} FS3D relies on structured Cartesian meshes and therefore decomposes the boundary of each cell into $4$ planar faces (i.e.\ edges) in two dimensions with normal vectors that are collinear with coordinate axes. The grid points used for the velocity components (face/edge centroids) are shifted with respect to the grid points used for the volume fractions (volume/area centroids) as shown in Figure~\ref{fig:staggered_grid}.\\

The idea of the operator splitting approach \cite{Rider1998,Rieber2004} is to decompose the full transport problem into a series of one-dimensional transport steps along the coordinate axis. Formally, this is achieved by decomposing the velocity as $v = (v_1,0) + (0,v_2) =: \tilde{v}_1 + \tilde{v}_2$. With this notation, equation \eqref{eqn:continuous_alpha_evolution_v2} reads as
\begin{equation}
\label{eqn:continuous_alpha_evolution_v3}
\begin{aligned}
\alpha_{ij}(t+\Delta t) = \left(\alpha_{ij}(t) - \frac{1}{|V_{ij}|} \int_{t}^{t+\Delta t} \int_{\partial V_{ij}} \chi \,  \tilde{v}_1 \cdot n \, dA \, d\tau + \frac{1}{|V_{ij}|} \int_{t}^{t+\Delta t} \int_{V_{ij}} \chi (\nabla \cdot \tilde{v}_1) \, dV \, d\tau \right) \\
- \frac{1}{|V_{ij}|} \int_{t}^{t+\Delta t} \int_{\partial V_{ij}} \chi \,  \tilde{v}_2 \cdot n \, dA \, d\tau + \frac{1}{|V_{ij}|} \int_{t}^{t+\Delta t} \int_{V_{ij}} \chi (\nabla \cdot \tilde{v}_2) \, dV \, d\tau.
\end{aligned}
\end{equation}
Note that the projected velocity fields $\tilde{v}_i$ are no longer divergence free and it is, therefore, important to consider the divergence correction terms in the discretization of \eqref{eqn:continuous_alpha_evolution_v3}. The discrete set of equation reads as
\begin{equation}
\label{eqn:alpha_update_formula}
\begin{aligned}
(1-\beta q_1^{i,j}) \alpha^*_{i,j} = &\alpha^n_{i,j} (1+(1-\beta)q_1^{i,j}) - \frac{\delta V_{i+1/2,j}-\delta V_{i-1/2,j}}{|V_{ij}|}, \\
(1-\beta q_2^{i,j}) \alpha^{n+1}_{i,j} = &\alpha^*_{i,j} (1+(1-\beta)q_2^{i,j}) - \frac{\delta V_{i,j+1/2}-\delta V_{i,j-1/2}}{|V_{ij}|},
\end{aligned}
\end{equation}
where $\delta V_{i \pm 1/2,j}$ denotes the volume flux over the edge $(i \pm 1/2,j)$ in a $x$-sweep, $\delta V_{i,j\pm 1/2}$ denotes the volume flux over the edge $(i,j \pm 1/2)$ in a $y$-sweep  and
\begin{align*}
q_1^{i,j} := \Delta t \, \frac{v_1(i,j) - v_1(i-1,j)}{\Delta x_1}, \quad q_2^{i,j} := \Delta t \, \frac{v_2(i,j+1) - v_2(i,j)}{\Delta x_2}.
\end{align*}
The choices $\beta=0$ and $\beta=1$ correspond to an explicit or implicit discretization of the divergence correction, respectively. For the present study, we choose $\beta=0.5$. The order of the direction of the sweeps is exchanged after each time step to avoid numerical asymmetries \cite{Strang1968}. After each directional split transport step, a heuristic volume redistribution algorithm similar to \cite{Harvie2000} is applied to enforce boundedness of the method, i.e. $0 \leq \alpha_{i,j}^n \leq 1$.

\paragraph{Approximation of the volume flux:} The volume fraction field $\alpha_{ij}$ is used to reconstruct a planar interface in each cell which is cut by the interface $\Sigma(t)$, in effect approximating $\chi$ with the Piecewise Linear Interface Calculation (PLIC). The reconstructed PLIC indicator $\tilde{\chi}$ is then used to approximate the flux integral
\[  \int_{t}^{t+\Delta t} \int_{A_f} \chi \,  \tilde{v}_d \cdot n \, dA \, d\tau, \quad d = 1,2, \]
where $A_f$ is a face (edge) of the cell $ij$. The geometrical interpretation of the above integral is a volume of the phase $\Omega^-$ (indicated with $\chi$) that passes through the face $A_f$. The calculation of this volume is performed by sweeping the face (edge) $A_f$ of the cell backwards along the trajectory given by $\tilde{v}_d$, and then clipping this swept volume with $\tilde{\chi}$. Note that the operator splitting simplifies the flux calculation because it simplifies the swept volume as a rectangular cuboid, whose intersection with $\tilde{\chi}$ is much simpler to calculate than for the swept volume that is traced along the Lagrangian trajectories given by the full velocity $v$. The simplification comes at the price of performing an additional PLIC interface intersection per splitting step.

\paragraph{Flux at boundary faces:} No special care is necessary to compute the flux over cell faces at an impermeable domain boundary since we assume a vanishing normal velocity leading to zero flux there.\\
\\
Accurate handling of the contact line motion by the geometrical VOF method, whose verification is outlined in section \ref{section:kinematics}, requires accurate interface reconstruction methods near the boundary of the solution domain. Proposed adaptations of well-known reconstruction algorithms (see Section~\ref{section:interface_reconstruction}) are described in detail in  Section~\ref{section:boundary_reconstruction}.

\subsection{Interface Reconstruction}
\label{section:interface_reconstruction}
A large number of methods have been developed for the geometrical reconstruction of the interface from the volume fraction field. An overview of reconstruction algorithms can be found in \cite{Pilliod2004,Aulisa2007}. In the present work, we consider two reconstruction algorithms, namely the classical method by Youngs \cite{Youngs1984} and the ELVIRA method due to Pilliod and Puckett \cite{Pilliod2004}. We propose extensions of these two methods to reconstruct an interface close to the boundary in Section~\ref{section:boundary_reconstruction}. To keep the formulas simple, assume that the mesh is equidistant in each direction with mesh sizes denoted as $\Delta x_1$ and $\Delta x_2$.\newline
\newline
A natural measure for the interface reconstruction error is the $L^1$-error defined as
\begin{align}
\mathcal{E}_1 := \int_\Omega |\chi(x) - \hat\chi(x) | \, dx,
\end{align}
where $\hat\chi$ is the characteristic function of the reconstructed domain.

\paragraph{Youngs Reconstruction Method:} The idea of the Youngs method is to approximate the interface normal vector by the discrete gradient of the volume fraction field, i.e.
\begin{align}
\label{eqn:youngs_normal}
 \nsigma \approx \nsigma^Y = - \frac{\nabla_h \alpha}{|\nabla_h \alpha|}.
\end{align}
Then a plane with orientation $\nsigma^Y$ is positioned such that the volume fraction in the local cell is matched (see \cite{Scardovelli2000} for details of the positioning algorithm). The gradient in \eqref{eqn:youngs_normal} is approximated by weighted central finite differences on a $3\times3$-block of cells. For an equidistant mesh, the gradient at cell $(i,j)$ is discretized with central finite differences as
\begin{equation}
\label{eqn:youngs_gradient}
\begin{aligned}
(\nabla_h \alpha)_1 = \ &\frac{1}{2} \frac{\alpha(i+1,j)-\alpha(i-1,j)}{2 \Delta x_1} + \frac{1}{4} \frac{\alpha(i+1,j+1)-\alpha(i-1,j+1)}{2 \Delta x_1} + \frac{1}{4} \frac{\alpha(i+1,j-1)-\alpha(i-1,j-1)}{2 \Delta x_1},\\
(\nabla_h \alpha)_2 = \ &\frac{1}{2} \frac{\alpha(i,j+1)-\alpha(i,j-1)}{2 \Delta x_2} + \frac{1}{4} \frac{\alpha(i+1,j+1)-\alpha(i+1,j-1)}{2 \Delta x_2} + \frac{1}{4} \frac{\alpha(i-1,j+1)-\alpha(i-1,j-1)}{2 \Delta x_2}.
\end{aligned}
\end{equation}
The Youngs Method is known to be one of the fastest methods for interface reconstruction from volume fractions fields. But it is only \emph{first-order} accurate with respect to the $L^1$-norm since it fails to reconstruct all planar interface exactly, see \cite{Pilliod2004}. This may be explained as a consequence of the lack of regularity of the volume fraction field $\alpha(i,j)$, which can be understood as the \emph{evaluation} at the cell centers of the continuous function
\[ F(x) = \frac{1}{|V_0|} \int_{V_0} \chi(x+x') \, dx', \]
obtained from averaging the phase indicator function over a control volume. It is well-known \cite{Scardovelli2000} that the latter function is only of class $\mathcal{C}^1$. One can, therefore, \emph{not} expect convergence of the finite differences scheme \eqref{eqn:youngs_gradient}. This underlines the need for more advanced methods.

\paragraph{(E)LVIRA Method:} The idea of the \textit{Least Squares VOF Interface Reconstruction Algorithm} (LVIRA) method proposed by Puckett \cite{Puckett.1991} is to find a planar interface reconstruction that \emph{minimizes} the quadratic deviations of the volume fractions in a $3\times3$-block under the constraint that this interface exactly reproduces the volume fraction in the central cell. Hence, one minimizes the functional
\begin{align}
\label{eqn:lvira_functional}
\mathcal{F} = \sum_{k,l=-1}^1[\tilde\alpha_{i+k,j+l}(n) - \alpha_{i,j}]^2,
\end{align}
where $\tilde\alpha_{i+k,j+l}(n)$ is the volume fraction in cell $(i+k,j+l)$ which is induced by a plane with orientation $n$ satisfying $\tilde\alpha_{ij}(n)=\alpha_{ij}$. Due to the nonlinear constraint, the minimization problem cannot be reformulated as a linear system of equations.\newline
\newline
Since the minimization of \eqref{eqn:lvira_functional} is computationally expensive compared to a direct method like the Youngs reconstruction, Pilliod and Puckett introduced the \textit{Efficient Least Squares VOF Interface Reconstruction Algorithm} (ELVIRA) \cite{Pilliod2004}. The computational costs are reduced by minimizing \eqref{eqn:lvira_functional} only over a \emph{finite} set of candidate orientations obtained in the following way:\\
\\
Suppose the interface can be described in the slope-intercept form
\begin{align} 
\label{eqn:slope_intercept_1}
x_2 = m_1 x_1 + b. 
\end{align}

Then the interface normal vector is either 
\begin{align*} 
\nsigma = (-m_1,1)/\sqrt{1+m_1^2} \quad \text{or} \quad \nsigma = (m_1,-1)/\sqrt{1+m_1^2}. 
\end{align*}
The slope is approximated by central-, forward- and backward-finite-differences of column sums. The candidates for the slope $m_x$ in the cell $(i,j)$ are
\begin{equation}
\label{eqn:elvira_column_sums}
\begin{aligned}
\pm m_1^b \quad \text{with} \quad m_1^b &= \frac{\Delta x_2}{\Delta x_1} \, \sum_{l=-1}^1 (\alpha_{i,j+l} - \alpha_{i-1,j+l}), \\
\pm m_1^c \quad \text{with} \quad m_1^c &= \frac{\Delta x_2}{2\Delta x_1} \, \sum_{l=-1}^1 (\alpha_{i+1,j+l} - \alpha_{i-1,j+l}), \\
\pm m_1^f \quad \text{with} \quad m_1^f &= \frac{\Delta x_2}{\Delta x_1} \, \sum_{l=-1}^1 (\alpha_{i+1,j+l} - \alpha_{i,j+l}).
\end{aligned}
\end{equation}
We observe that the slope should be approximated with $+m_1^b$, $+m_1^c$ or $+m_1^f$ if the second component of $\nsigma$ is positive and vice versa with $-m_1^b$, $-m_1^c$ or $-m_1^f$ if the second component of $\nsigma$ is negative. This results in $6$ candidates for the normal vector obtained from column sums in the $x_2$-direction.\\
\\
Obviously, it is not always possible to represent the interface as a graph over $x_1$. Therefore, one also has to consider the case
\begin{align} 
\label{eqn:slope_intercept_2}
x_1 = m_2 x_2 + b. 
\end{align}
This gives rise to analogous approximations for $m_2^b$,$m_2^c$ and $m_2^f$. This results in  $12$ candidates for the interface normal in two dimensions. It can be shown to be sufficient to reconstruct any straight line exactly which makes the method formally second-order accurate with respect to the $L^1$-error, see \cite{Pilliod2004}.

\subsection{Interface Reconstruction close to the Boundary}
\label{section:boundary_reconstruction}
We propose an adaptation of well-known methods for interface reconstruction in Volume-of-Fluid methods, aiming for accurate reconstruction of the interface close to the domain boundary. The FORTRAN implementations of the developed interface reconstruction methods are available online in an open research data repository, see \cite{OnlineRepoVOF}.\\
\\
For simplicity of notation, we consider in the following the domain boundary at $x_2 = 0$, i.e.\ we consider an interface cell with index $(i,j=1)$.

\subsubsection{Boundary Youngs Method}
\begin{figure}[ht]
 \centering
 \subfigure{\includegraphics[width=4cm]{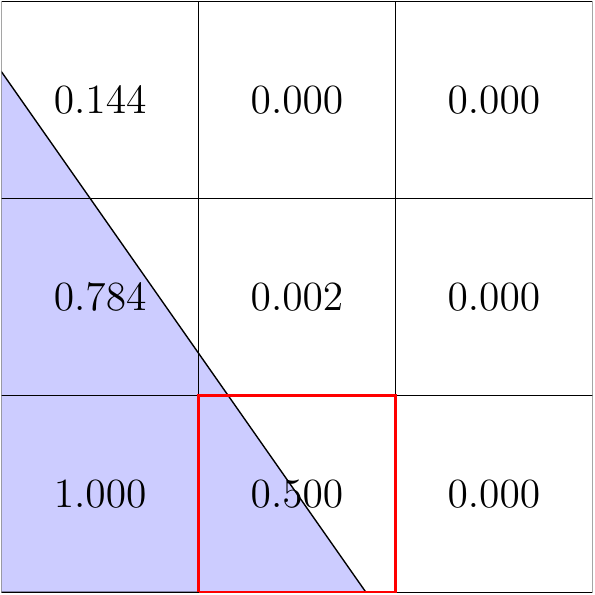}}
 \subfigure{\includegraphics[width=4cm]{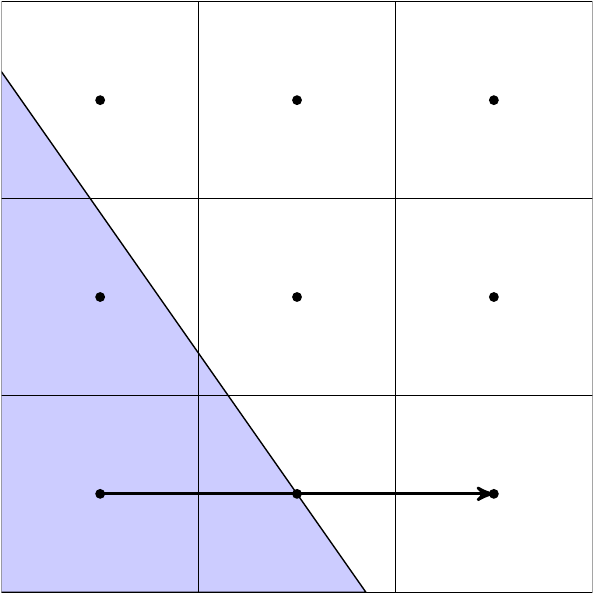}}
 \subfigure{\includegraphics[width=4cm]{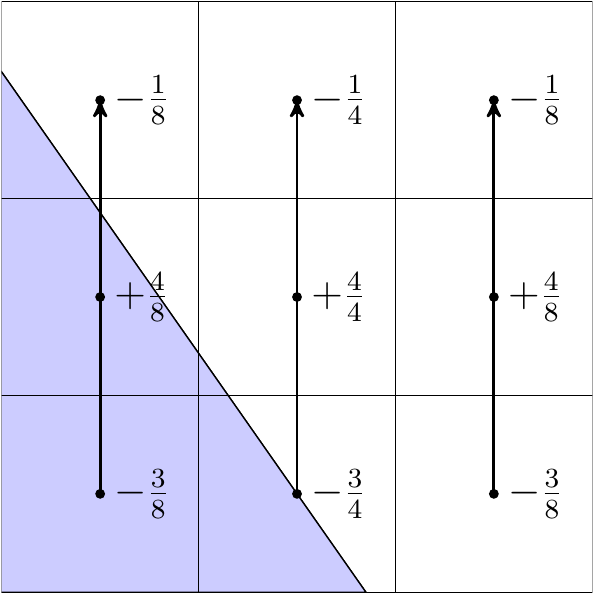}}
 \subfigure{\includegraphics[width=4cm]{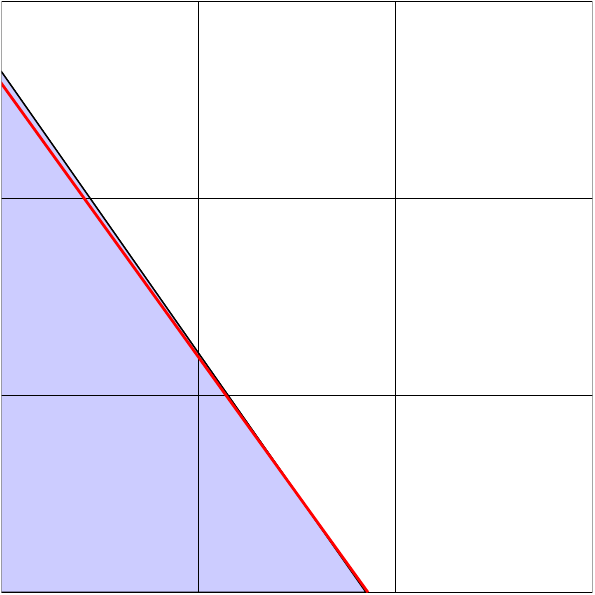}}
 \caption{\emph{Boundary} Youngs Reconstruction Method for an equidistant mesh.}
 \label{fig:boundary_youngs}
\end{figure}

We consider a $3\times3$-block of cells and aim at reconstructing the interface in the lower middle cell. We propose to discretize the gradient of the volume fraction field in two space dimensions in the following way (see Figure~\ref{fig:boundary_youngs}): 
\begin{enumerate}[(a)]
 \item Tangential to the domain boundary central finite differences are used.
 \item In normal direction to the domain boundary, weighted \emph{forward} finite differences are employed.
\end{enumerate}
From Taylor's formula, one can show that for a $\mathcal{C}^3$-function $f$, the first derivative $f'(x)$ can be approximated with second-order accuracy according to
\begin{align*}
f'(x) = \frac{-f(x+2\Delta x) + 4 f(x+\Delta x) - 3 f(x)}{2 \Delta x} + \mathcal{O}(\Delta x^2).
\end{align*}
This formula is applied to approximate the derivative of the volume fraction normal to the boundary. Note, however, that the volume fraction $\alpha$ is only $\mathcal{C}^1$ (see Section~\ref{section:interface_reconstruction}). One can, therefore, \emph{not} expect a convergence of the orientation with that method. But we still consider it here since it is a straightforward extension of the widely used Youngs method to the boundary case.\newline
\newline
For an equidistant grid in two space dimensions, the Boundary Youngs gradient in a cell with index $(i,1)$ is discretized as
\begin{equation}
\label{eqn:boundary_youngs_gradient}
\begin{aligned}
(\nabla_h \alpha)_1 = \ &\frac{\alpha(i+1,1)-\alpha(i-1,1)}{2\Delta x_1},\\
(\nabla_h \alpha)_2 = \ & \frac{-\alpha(i,3)+4 \alpha(i,2)-3 \alpha(i,1)}{4 \, \Delta x_2}\\
+ \ & \frac{-\alpha(i+1,3)+4\alpha(i+1,2)-3\alpha(i+1,1)}{8 \, \Delta x_2}\\
+ \ & \frac{-\alpha(i-1,3)+4\alpha(i-1,2)-3\alpha(i-1,1)}{8 \, \Delta x_2}.
\end{aligned}
\end{equation}

\subsubsection{Boundary ELVIRA Method}
\begin{figure}[ht]
 \centering
 \includegraphics[width=0.5\columnwidth]{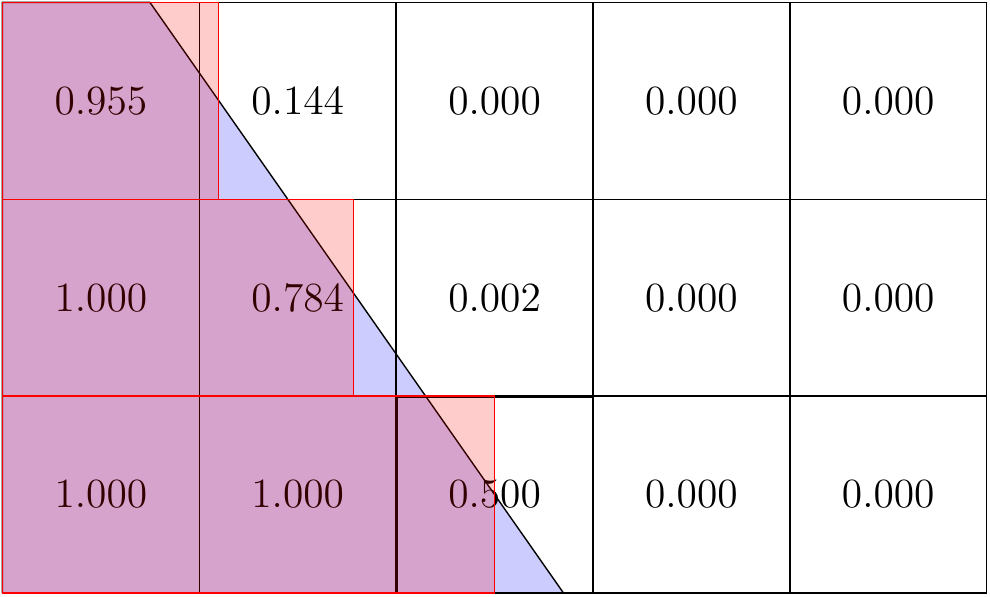}
 \caption{Boundary ELVIRA method on a $5\times3$-stencil in 2D.}
 \label{fig:boundary_elvira}
\end{figure}

In order to allow for mesh convergent results for the contact angle evolution, one needs a reconstruction method which is second-order accurate at the boundary. Therefore, we propose the following adaptation of the ELVIRA method due to Pilliod and Puckett \cite{Pilliod2004}: Minimize the functional
\begin{align}
\label{eqn:boundary_lvira_functional}
\mathcal{F}_b = \sum_{k=-2}^2 \sum_{l=0}^2[\tilde\alpha_{i+k,1+l}(n) - \alpha_{i,1}]^2,
\end{align}
where $\tilde\alpha_{i+k,1+l}(n)$ is the volume fraction in cell $(i+k,1+l)$ which is induced by a plane with orientation $n$ satisfying $\tilde\alpha_{i,1}(n)=\alpha_{i,1}$. Here the minimization is performed over a \emph{larger} stencil of $5\times3$ cells. This turns out to be necessary to reconstruct every straight line at the boundary exactly. Following the idea of the \textit{Efficient Least Squares VOF Interface Reconstruction Algorithm} \cite{Pilliod2004}, the functional \eqref{eqn:boundary_lvira_functional} is minimized over a finite set of candidate orientations obtained from finite differences of column sums (see Figure~\ref{fig:boundary_elvira}, where the column sums in horizontal direction are visualized in red).\\

The following candidate slopes are computed in the cell $(i,1)$ in normal direction to the boundary
\begin{equation}
\label{eqn:boundary_elvira_column_sums_normal}
\begin{aligned}
m_1^c &= \frac{\Delta x_2}{2\Delta x_1} \, \sum_{l=0}^2 (\alpha_{i+1,1+l} - \alpha_{i-1,1+l}), \\
m_1^b &= \frac{\Delta x_2}{\Delta x_1} \, \sum_{l=0}^2 (\alpha_{i,1+l} - \alpha_{i-1,1+l}), \\
m_1^{b^\ast} &= \frac{\Delta x_2}{\Delta x_1} \, \sum_{l=0}^2 (\alpha_{i-1,1+l} - \alpha_{i-2,1+l}), \\
m_1^f &= \frac{\Delta x_2}{\Delta x_1} \, \sum_{l=0}^2 (\alpha_{i+1,1+l} - \alpha_{i,1+l}),\\
m_1^{f^\ast} &=  \frac{\Delta x_2}{\Delta x_1} \, \sum_{l=0}^2 (\alpha_{i+2,1+l} - \alpha_{i+1,1+l}).
\end{aligned}
\end{equation}
The following candidate slope is computed from sums tangentially to the boundary
\begin{equation}
\label{eqn:boundary_elvira_column_sums_tangential}
\begin{aligned}
m_2^f &= \frac{\Delta x_1}{\Delta x_2} \, \sum_{l=-2}^2 (\alpha_{i+l,2}-\alpha_{i+l,1}), \\
m_2^{f^\ast} &= \frac{\Delta x_1}{2\Delta x_2} \, \sum_{l=-2}^2 (\alpha_{i+l,3}-\alpha_{i+l,1}),\\
m_2^{f^{\ast\ast}} &= \frac{\Delta x_1}{\Delta x_2} \, \sum_{l=-2}^2 (\alpha_{i+l,3}-\alpha_{i+l,2}).
\end{aligned}
\end{equation}
This yields $16$ candidates for the interface normal. We can demonstrate by numerical experiments that this is sufficient to reconstruct any straight line at the boundary up to machine precision, see \cite{OnlineRepoVOF}.

\subsubsection{Numerical errors} 
\begin{figure}[ht]
 \centering
 \subfigure[Translation test away from the boundary.]{\includegraphics[width=0.45\columnwidth]{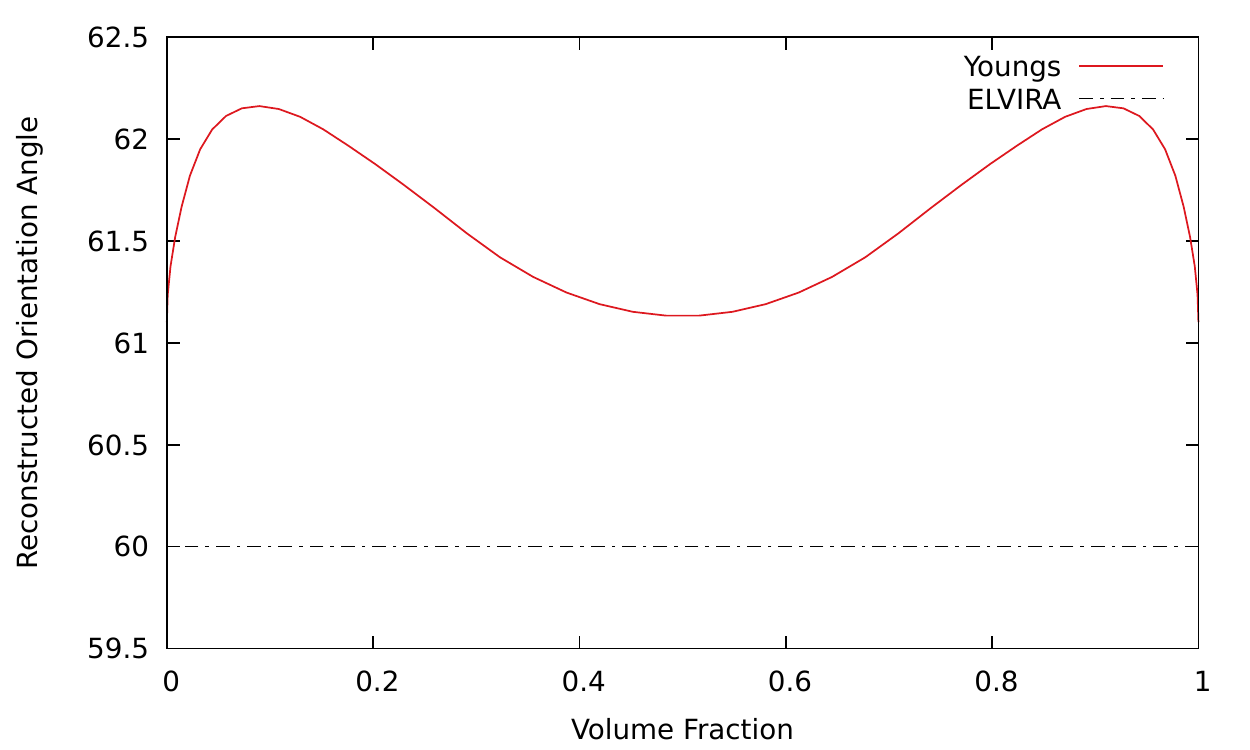}\label{fig:translation_test_bulk}}
 \subfigure[Translation test at the boundary.]{\includegraphics[width=0.45\columnwidth]{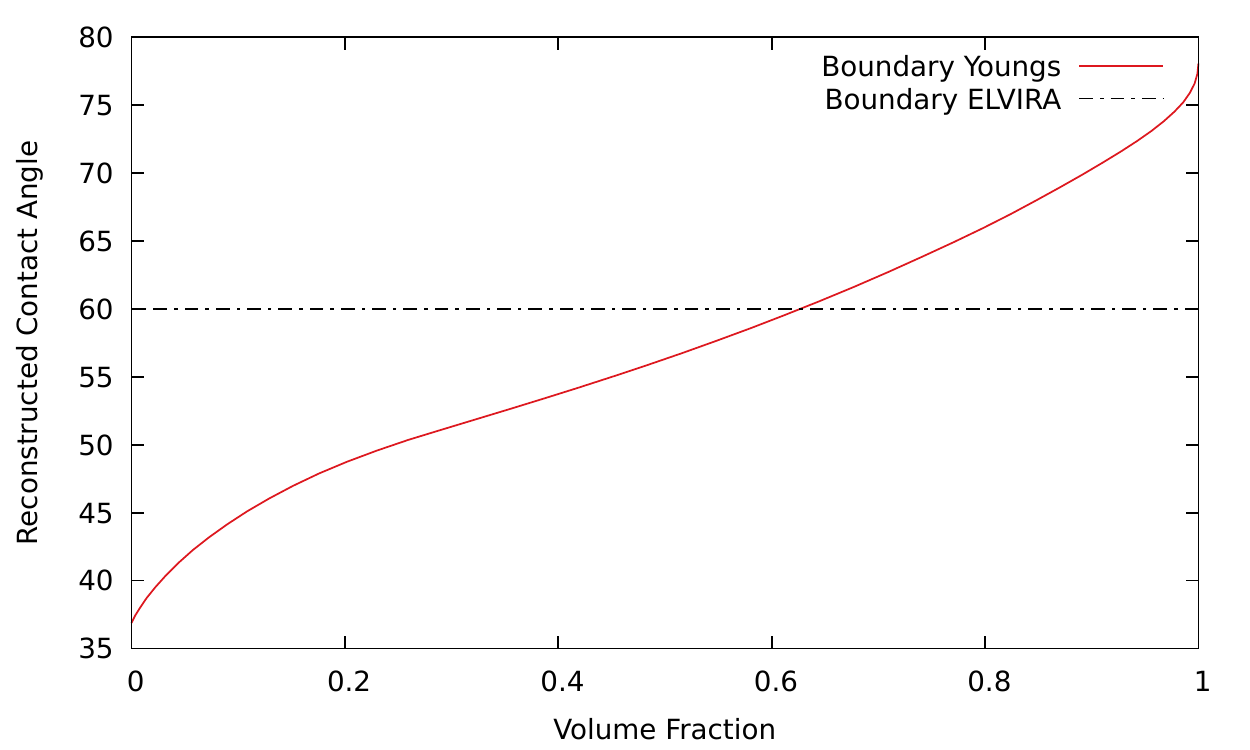}\label{fig:translation_test_boundary}}
 \caption{Translation of a straight line with a fixed orientation angle of $\theta=60^\circ$.}
 \label{fig:translation_test}
\end{figure}

It is well-known from the literature that the standard Youngs method \eqref{eqn:youngs_gradient} fails to reconstruct arbitrary straight lines, while the error is typically of the order of a few degrees. As a numerical test, a straight line is moved with a fixed inclination angle on an equidistant grid with $\Delta x_1 = \Delta x_2$. This motion produces volume fractions ranging from $0$ to $1$ for the considered computational cell away from the boundary. The reconstructed orientation with the standard Youngs and ELVIRA methods are shown in Figure~\ref{fig:translation_test_bulk}. While the ELVIRA method always delivers the correct angle, the Youngs method shows an error of about $1^\circ-2^\circ$ in the considered example. The situation is much different for the same translation test for a \emph{boundary} cell, see Figure~\ref{fig:translation_test_boundary}. While the Boundary ELVIRA method is still able to deliver the correct orientation, the Boundary Youngs method shows a large error of up to $\pm 20^\circ$ that is also highly dependent on the position of the interface. Therefore, one can only expect a very rough estimate of the contact angle from the Boundary Youngs method which cannot converge with mesh refinement.

\section{Results}
\begin{figure}[ht]
 \centering
 \subfigure[Regular case.]{\includegraphics[width=0.35\columnwidth]{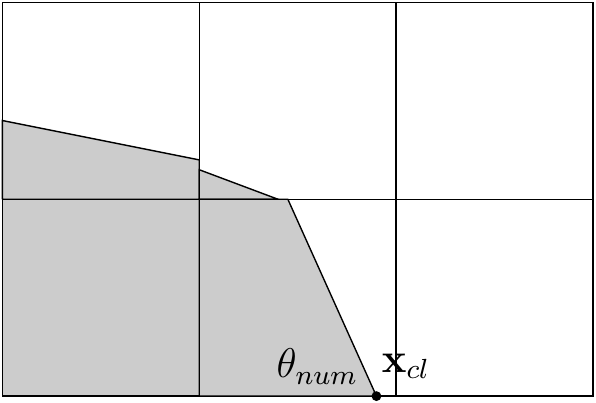}}
 \subfigure[Irregular case, no contact point detected.]{\includegraphics[width=0.35\columnwidth]{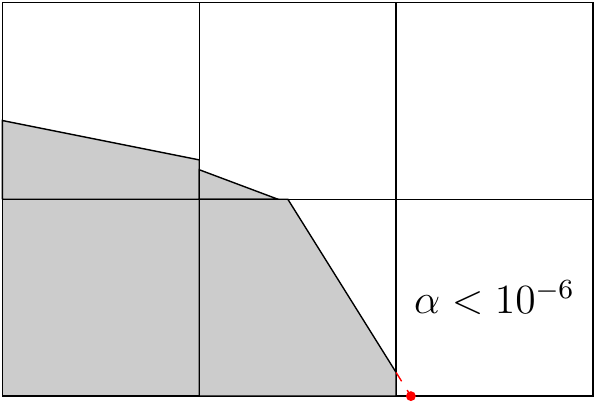}\label{fig:data_from_plic_execption}}
 \caption{Contact angle and contact line position from PLIC reconstruction.}
 \label{fig:data_from_plic}
\end{figure}

To verify the advective transport of the contact angle, a spherical cap sitting at the boundary is initialized and transported using different velocity fields. To study the convergence in space and time, the time step $\Delta t$ is linked to the grid spacing $\Delta x$ by fixing the Courant number
\[ \cfl = \frac{ \Delta t \, \norm{v}_{L^\infty(\Omega)}}{\Delta x}. \]
The influence of the choice of the Courant number is discussed below (see Figures~\ref{fig:cfl_study_1} and \ref{fig:cfl_study_2}). The following computational examples are carried out with (unless stated otherwise)
\[ \cfl = 0.2. \]
Fixing the Courant number defines a temporal grid $\mathcal{T}$ (equidistant if $v$ does not depend on time). We report the error for both the contact line position and the contact angle in the maximum norm over all time steps, i.e.
\begin{align*} 
\errortheta([0,T]) := \max_{t_i \in \mathcal{T}\cap [0,T]} |\thetanum(t_i) - \thetaref(t_i) | \quad \text{and} \quad \errorcl([0,T]) := \max_{t_i \in \mathcal{T}\cap [0,T]} |(\xclnum(t_i) - \xclref(t_i))/R_0|,
\end{align*}
as a function of $\Delta x/R_0$. Note that the error in the contact line position is normalized by the initial radius $R_0$. The reference values $\xclref(t)$ and $\thetaref(t)$ come either from an exact or from a numerical solution of the ordinary differential equations \eqref{eqn:characteristics_ode} and \eqref{eqn:theta_evolution_equation}. The numerical values for the contact line position $\xclnum$ and the contact angle $\thetanum$ are evaluated directly from the reconstructed PLIC element intersecting the domain boundary (see Figure~\ref{fig:data_from_plic}). To this end, the point of intersection of the local interface with the domain boundary is computed. If this point lies within the cell, the cell is recognized as a contact line cell and the contact angle and the contact line position are computed. Note that, due to the finite reconstruction tolerance of the VOF method (in this case $10^{-6}$) , irregular cases where no contact point is found may occur. An example is sketched in Figure~\ref{fig:data_from_plic_execption}, where the point of intersection lies slightly outside the current cell but the volume fraction of the neighbor cell is below the reconstruction tolerance so that it is not recognized as an interface cell. These irregular cases are excluded from the following error analysis.

\begin{figure}[ht]
 \centering
 \includegraphics[width=1.0\columnwidth]{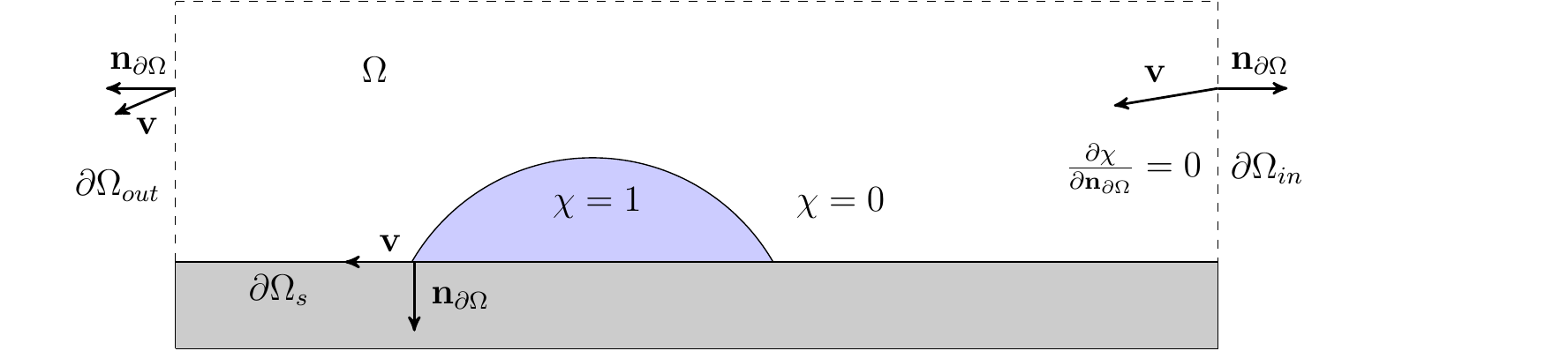}
 \caption{Computational setup.}
 \label{fig:computational_setup}
\end{figure}

\paragraph{Computational setup:} For the subsequent examples, we choose the following common setup. The computational domain $\Omega = [0,1]\times[0,0.25]$ is covered by an equidistant Cartesian grid of $N \times N/4$ cells, where $N$ varies from $128$ to $2048$. A spherical cap with dimensionless radius $R_0 = 0.2$ is initialized at the ``solid boundary''
\[ \partial\Omega_s = \{ (x_1,0): 0 \leq x_1 \leq 1 \}, \]
where the flow is assumed to be tangential. The center of the sphere is placed at $(0.4,-0.1)$ yielding an initial contact angle of 
\[ \theta_0 = \arccos\left(\frac{0.1}{0.2}\right) = \frac{\pi}{3}. \]
Since we are only interested in the \emph{local} transport of the interface, we can allow for an artificial inflow boundary to the computational domain away from the contact line. Inflow boundaries are characterized by $v \cdot \ndomega < 0$ (see Figure~\ref{fig:computational_setup}). At inflow boundaries we formally apply a homogeneous Neumann boundary condition for the phase indicator function, i.e.
\begin{align} 
\label{eqn:inflow_boundary}
\frac{\partial \chi}{\partial \ndomega} = 0 \quad  \text{on} \quad \partial\Omega_{in}(t) = \{x \in \partial\Omega: \, v(t,x) \cdot \ndomega(x) < 0\}. 
\end{align}
This condition is straightforward to implement using a simple constant continuation of the volume fraction field into a layer of ``ghost cells'' (see, e.g., \cite{LeVeque2002}). Here we only consider the case where the interface does not meet the artificial boundary such that \eqref{eqn:inflow_boundary} simply states that no additional volume is transported into the computational domain. In particular, the boundary condition does not affect the dynamics of the interface and \eqref{eqn:theta_evolution_equation} still holds.\\

Three examples for the transporting velocity field are studied with the Youngs and ELVIRA methods, where these methods are combined with their newly developed boundary versions to treat the contact line advection.

\subsection{Vortex-in-a-box test}
\begin{figure}[ht]
 \centering
 \includegraphics[width=0.6\columnwidth]{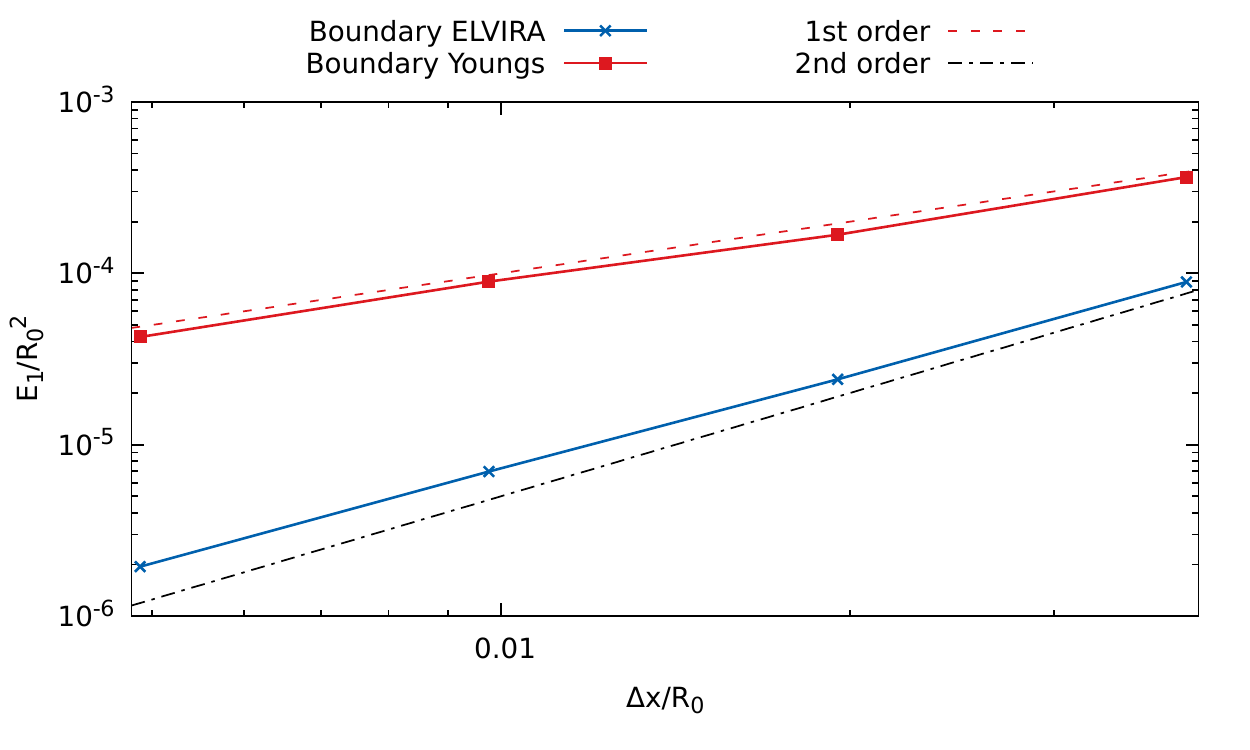}
 \caption{Convergence with respect to the discrete $L^1$-norm for the field \eqref{eqn:2D_shear_field} comparing initial and final shapes.}
 \label{fig:l1_error}
\end{figure}

We start with a classical test for interface advection methods given by
\begin{align}
\label{eqn:2D_shear_field}
v(x_1,x_2) = v_0 \cos\left( \frac{\pi t}{\tau} \right) (-\sin(\pi x_1)\cos(\pi x_2), \cos(\pi x_1)\sin(\pi x_2)).
\end{align}
This particular field called ``vortex-in-a-box'' has been routinely used to test numerical methods for interface advection; see \cite{Rider1998}, \cite{Tryggvason2011}. In the classical test this velocity field is used to strongly deform a sphere into a spiral. Due to the periodicity in time, it follows that the initial shape at $t=0$ and the final shape at $t = \tau$ would coincide if the problem is solved exactly. This allows to study aspects of the convergence behavior of the advection method even though the solution to the advection problem with the velocity field \eqref{eqn:2D_shear_field} is not known. The discrete $L^1$-error
\begin{align} 
\label{eqn:l1_error}
E_1 = \sum_{ij} |\alpha_{ij}(\tau) - \alpha_{ij}(0)| \, \Delta x_1 \Delta x_2 
\end{align}
is usually used to quantify the rate-of-convergence. Note, however, that this kind of test does not say anything about the intermediate \emph{dynamics} of the numerical solution\footnote{Formally, a numerical method which keeps the volume fractions fixed passes this test with zero error.}. Here we revisit this classical test in the presence of a moving contact line. The results for $v_0 = 0.1$ and $\tau = 0.2$ are reported in Figure~\ref{fig:l1_error}. The simulations are carried out with a fixed Courant number of $\cfl = 0.2$, where the numerical time step is chosen such that $t=\tau$ is reached after an integer number of time steps. As expected from the case without a contact line, the Boundary Youngs method shows a first-order convergence while the Boundary ELVIRA method is nearly second-order convergent.\\
\\
Thanks to the kinematic evolution equation, it is also possible to study the \emph{dynamics} of the advection in terms of the contact line position and the contact angle. The ordinary differential equations \eqref{eqn:characteristics_ode} and \eqref{eqn:theta_evolution_equation} are solved numerically to obtain reference solutions $\xclref(t)$ and $\thetaref(t)$.

\begin{figure}[ht]
  \centering
  \subfigure{\includegraphics[width=0.45\columnwidth]{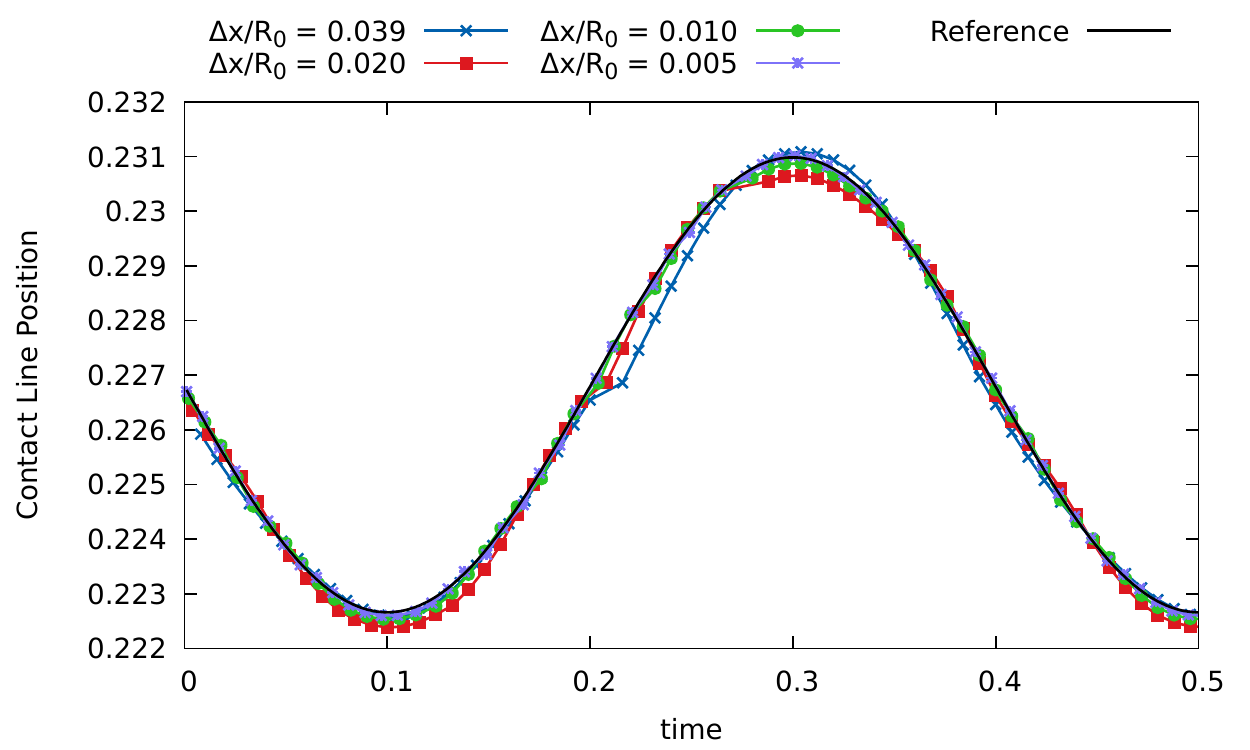}}
  \subfigure{\includegraphics[width=0.45\columnwidth]{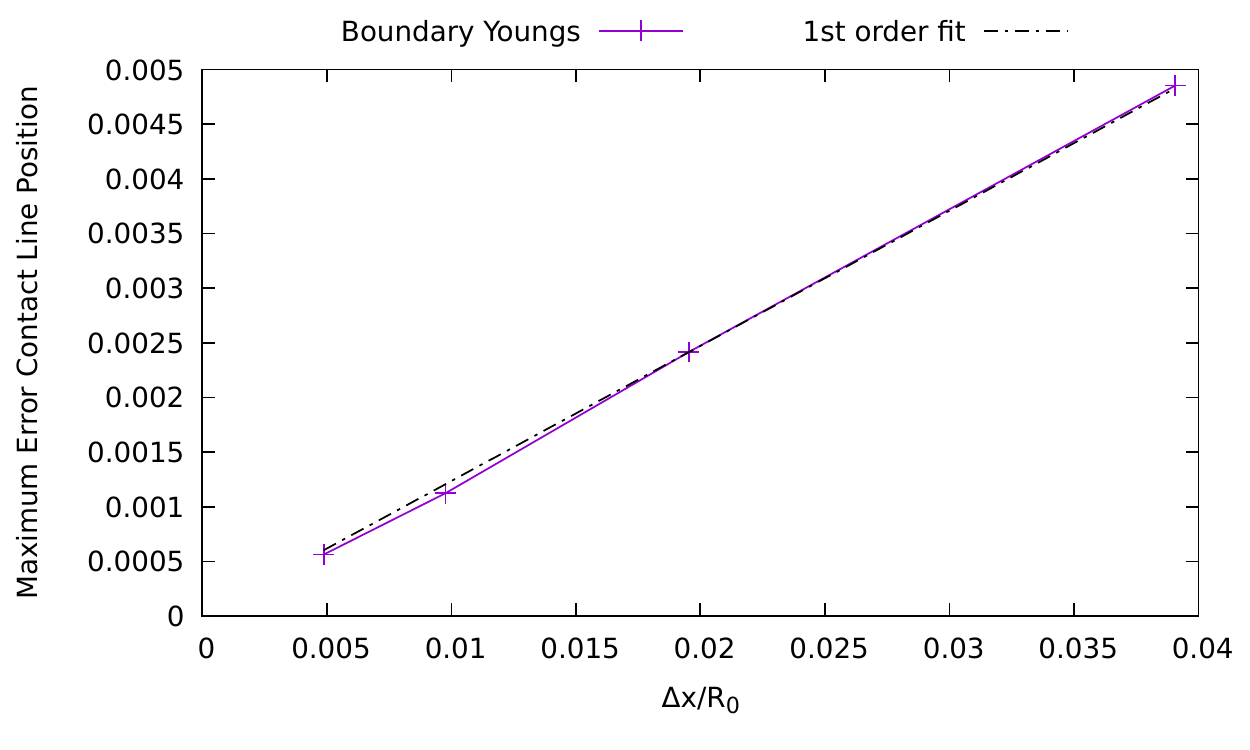}}
  \caption{Numerical motion of the contact line for the field \eqref{eqn:2D_shear_field} using the Boundary Youngs reconstruction.}
  \label{fig:plot_shear_youngs_position_2d}
 \end{figure}

\begin{figure}[ht]
  \centering
  \subfigure{\includegraphics[width=0.45\columnwidth]{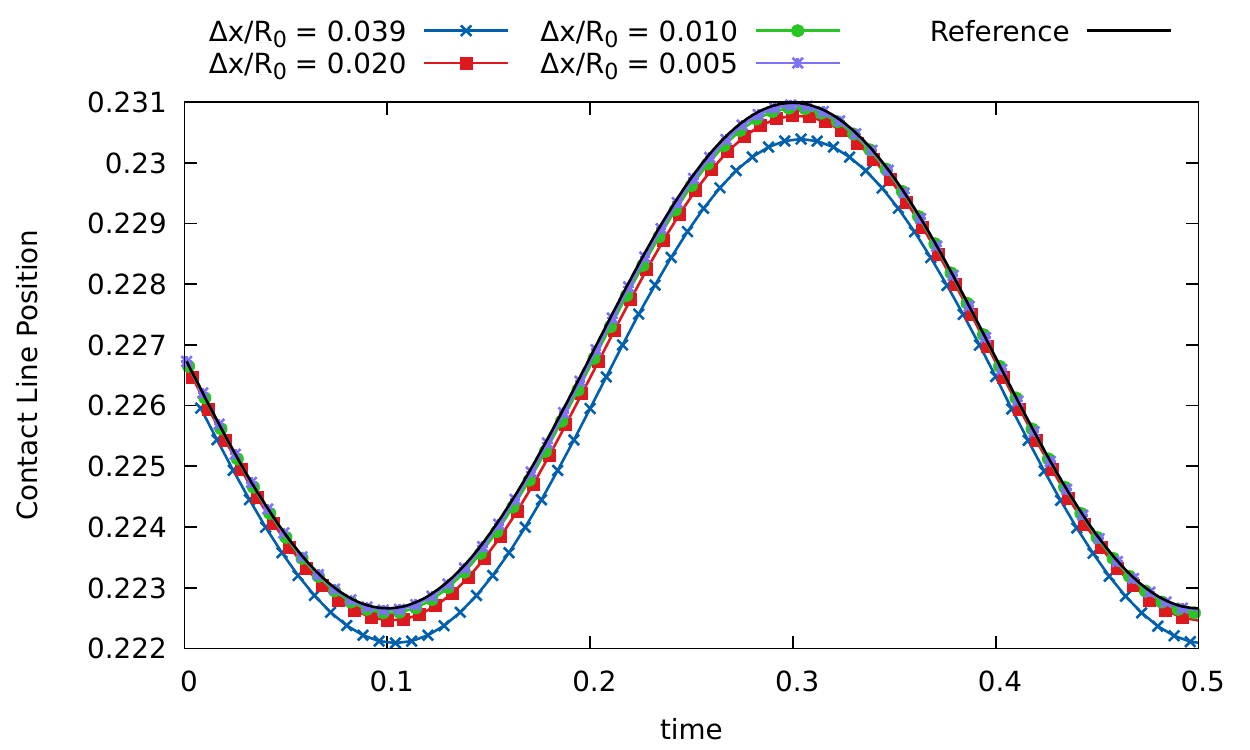}}
  \subfigure{\includegraphics[width=0.45\columnwidth]{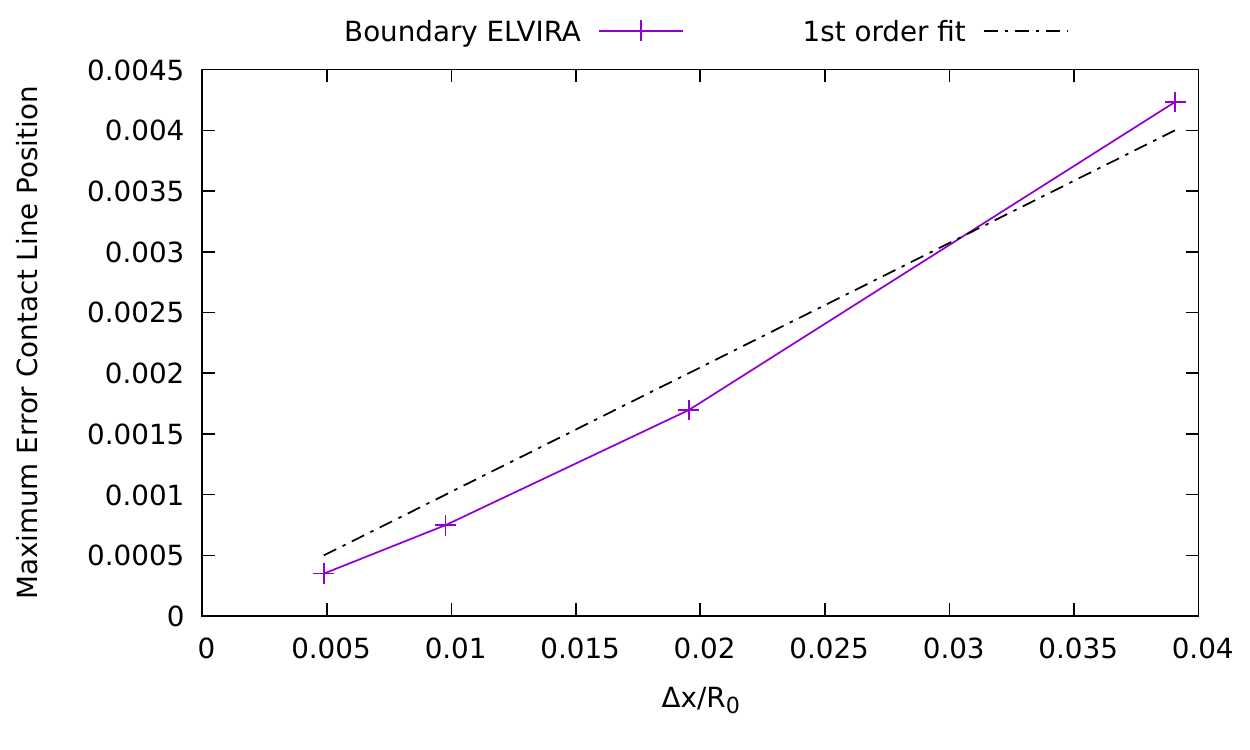}}
  \caption{Numerical motion of the contact line for the field \eqref{eqn:2D_shear_field} using the Boundary ELVIRA reconstruction.}
  \label{fig:plot_shear_elvira_position_2d}
 \end{figure}

\paragraph{Contact Line Motion:} The numerical evolution of the left (in this case the advancing) contact point reconstructed from the PLIC interface is investigated for the Boundary Youngs and Boundary ELVIRA method. It is found that both the Boundary Youngs and Boundary ELVIRA method deliver at least first-order convergent results for the motion of the contact line, see Figures~\ref{fig:plot_shear_youngs_position_2d} and \ref{fig:plot_shear_elvira_position_2d}.

\begin{figure}[ht]
  \centering
  \subfigure{\includegraphics[width=0.45\columnwidth]{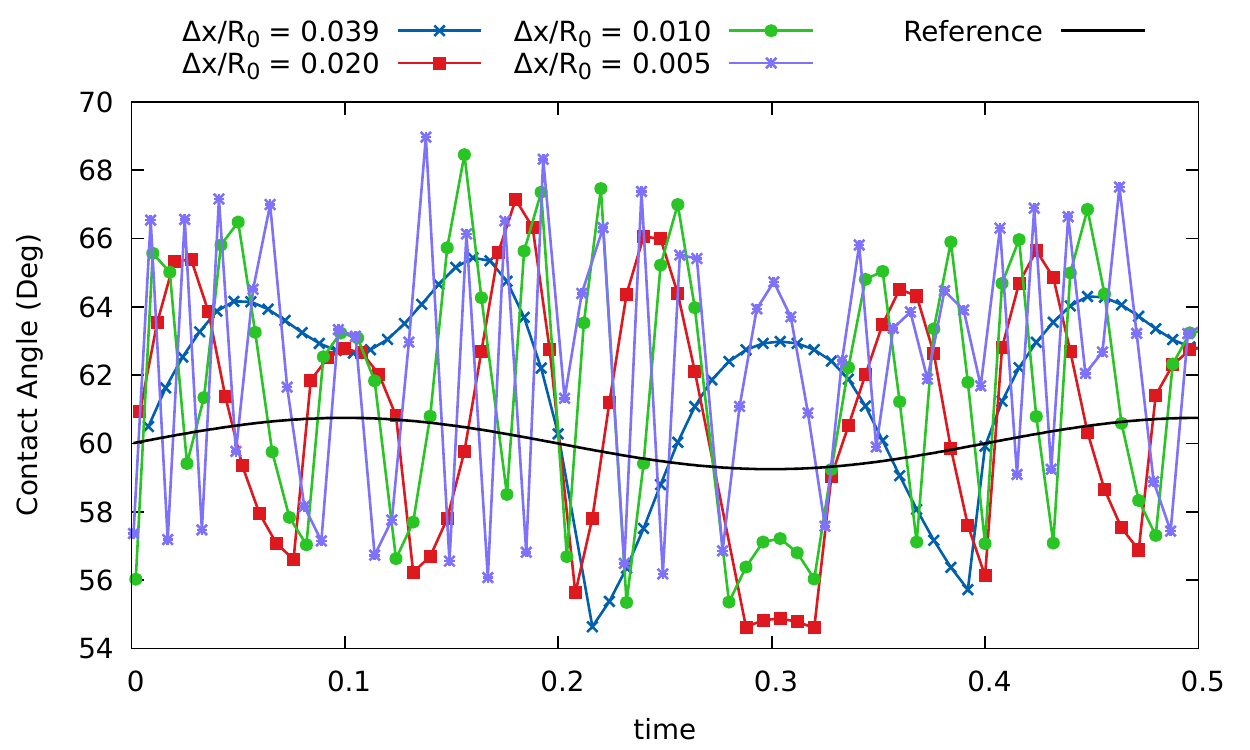}}
  \subfigure{\includegraphics[width=0.45\columnwidth]{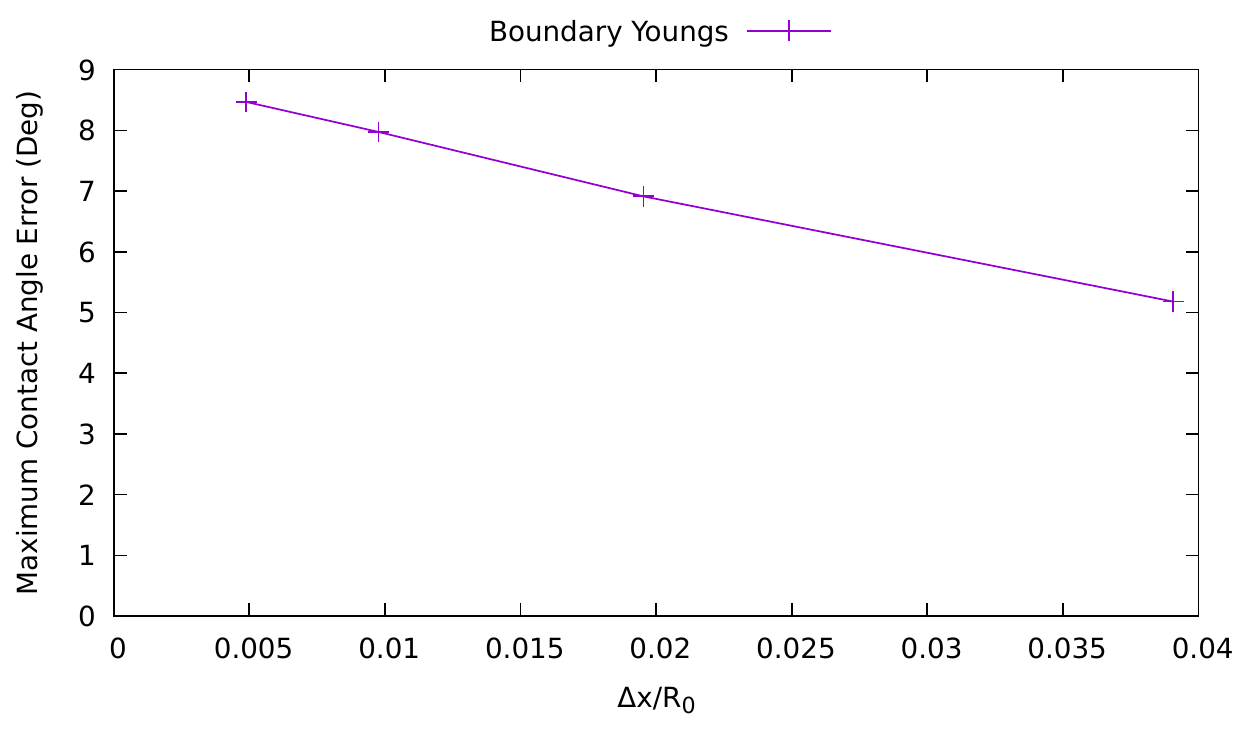}}
  \caption{Numerical $\theta$ evolution for the field \eqref{eqn:2D_shear_field} using the Boundary Youngs reconstruction.}
  \label{fig:results_shear_youngs_2d}
 \end{figure}

\begin{figure}[ht]
  \centering
 \subfigure{\includegraphics[width=0.45\columnwidth]{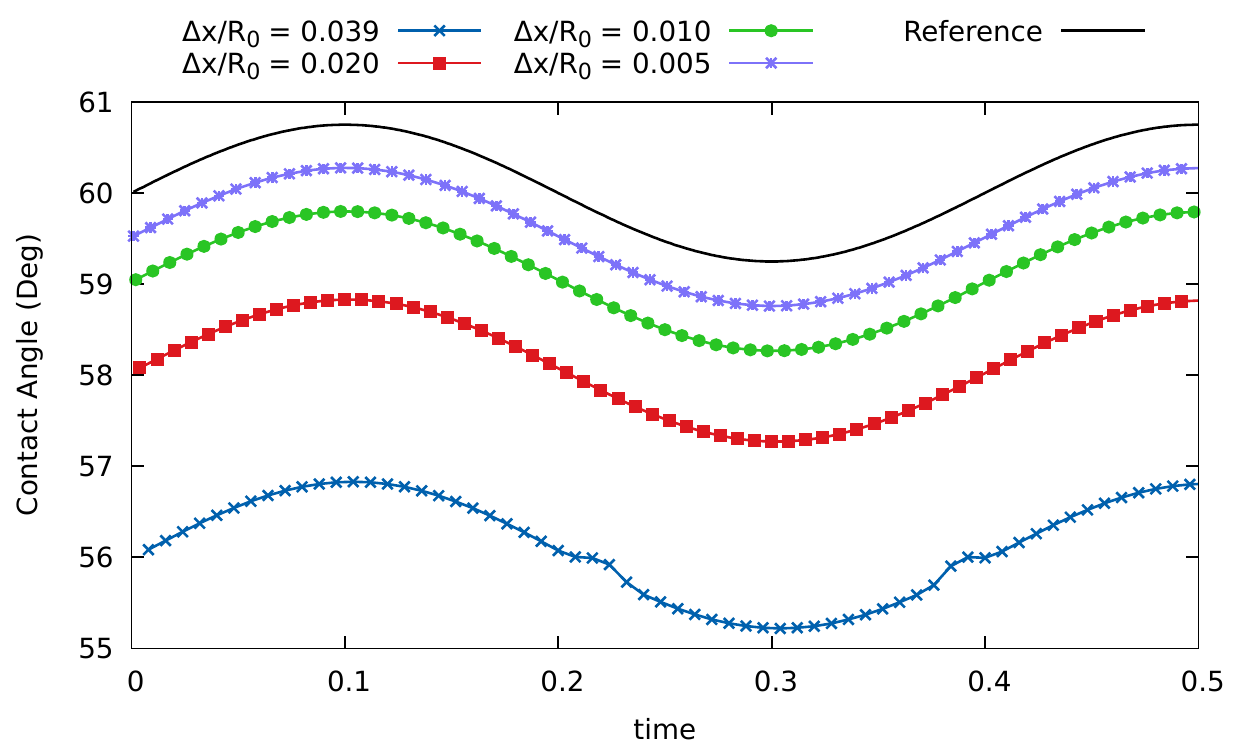}}
 \subfigure{\includegraphics[width=0.45\columnwidth]{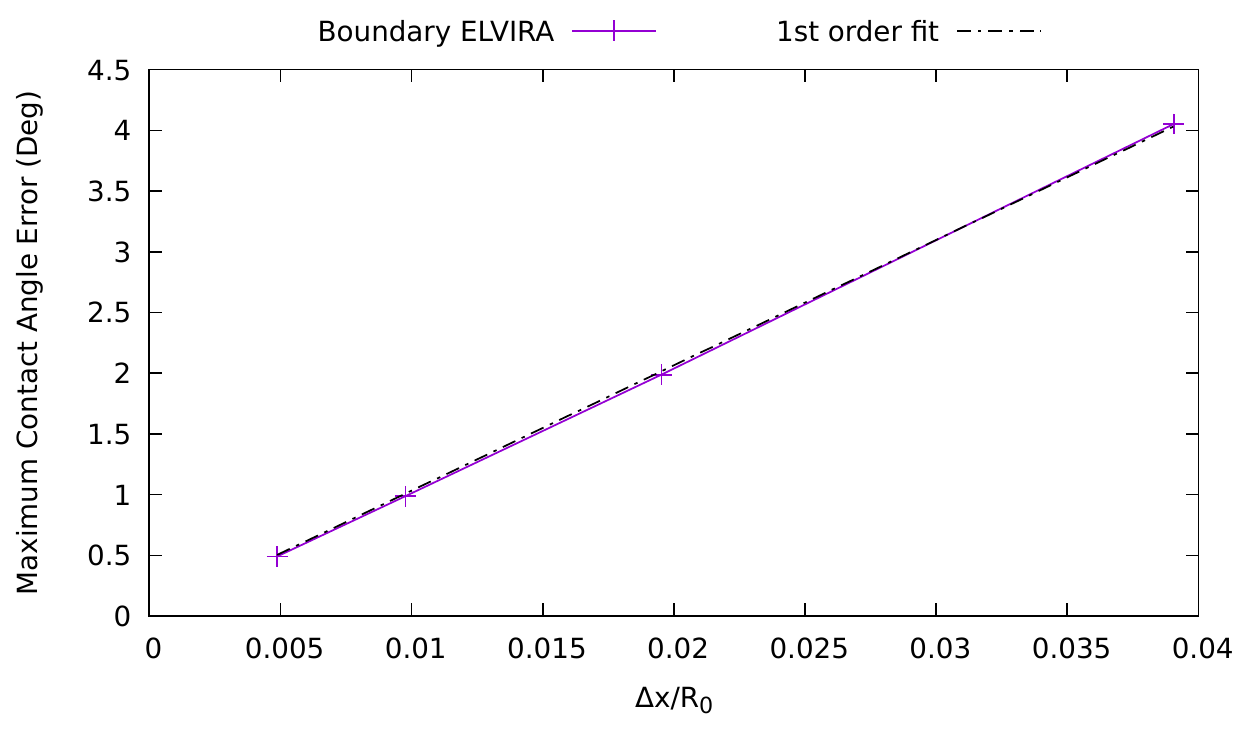}}
 \caption{Numerical $\theta$ evolution for the field \eqref{eqn:2D_shear_field} using the Boundary ELVIRA method.}
 \label{fig:results_shear_elvira_2d}
 \end{figure}

\begin{figure}[ht]
 \centering
 \includegraphics[width=0.5\columnwidth]{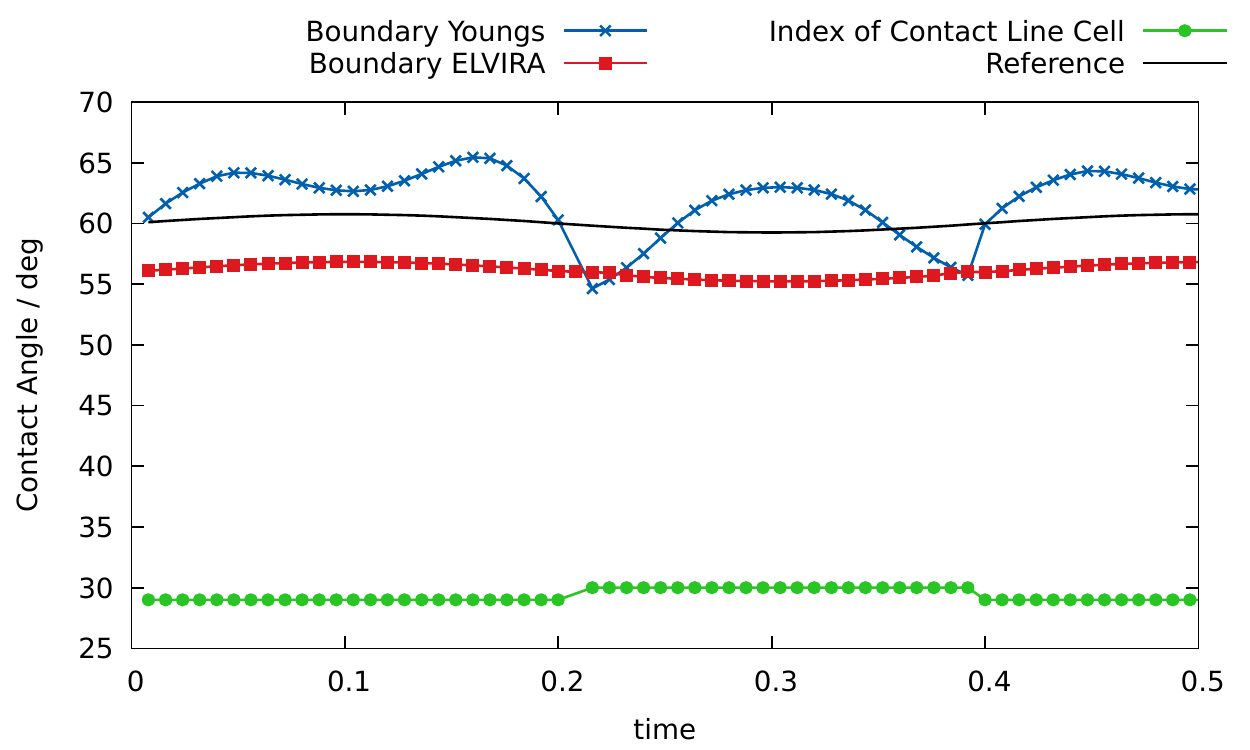}
 \caption{Jump in $\thetanum$ for the field \eqref{eqn:2D_shear_field} (here for $\Delta x/R_0 = 0.039$).}
 \label{fig:contact_angle_jump}
\end{figure}

\paragraph{Contact Angle Evolution:} The results for the contact angle over time are shown in Figures~\ref{fig:results_shear_youngs_2d} and \ref{fig:results_shear_elvira_2d}. While the numerical solution convergences to the reference solution for the Boundary ELVIRA method, the Boundary Youngs method does, as expected, \emph{not} deliver mesh convergent results. In fact, one observes a strong oscillation of the reconstructed contact angle with a jump discontinuity when the contact line passes from one cell to the other; see Figure~\ref{fig:contact_angle_jump}, where the reconstructed contact angle is plotted along with the discrete cell index. This behavior might be due to the spatial structure of the reconstruction error as reported in Figure~\ref{fig:translation_test}. Clearly, the frequency of these jumps increases with mesh refinement leading to the strongly oscillatory behavior. The error in the maximum norm may even increase with mesh refinement, see Figure~\ref{fig:results_shear_youngs_2d}. Therefore, the Boundary Youngs method does \emph{not} allow for a meaningful evaluation of the contact angle based on the local interface orientation even though it is first-order convergent with respect to the contact line motion and the discrete $L^1$-error regarding the initial and final shape comparison.\\
\\
Following Figures~\ref{fig:results_shear_elvira_2d} and \ref{fig:contact_angle_jump}, the evolution of the numerical contact angle for the Boundary ELVIRA method is reasonably smooth even on coarse grids. Some small oscillations are visible which, however, disappear with mesh refinement. In fact, the method shows a \emph{first-order} convergence with respect to $\errortheta([0,0.5])$. The maximum error on the finest mesh with $\Delta x/R_0 = 5 \cdot 10^{-3}$ is about $0.5$ degrees.

\subsection{Linear Velocity Field}
We now consider linear velocity fields of the form \eqref{eqn:general_linear_field}. In this case, the explicit solution \eqref{eqn:exact_solution_2d} is available for verification. We choose the example
\begin{align}
\label{eqn:2D_navier_field}
v(x_1,x_2) = \left(-0.2 + 0.1 x_1 - 2 x_2, \, -0.1 x_2 \right).
\end{align}
The time evolution is investigated up to dimensionless time $T=0.4$. According to \eqref{eqn:exact_solution_2d}, the exact solution for the left contact point is given by
\begin{align*}
x_1(t) &= x_1^0 e^{0.1t} - 2(e^{0.1t} -1), \quad \thetaref(t) = \frac{\pi}{2} + \arctan\left(-\frac{1}{\sqrt{3}} e^{0.2 t} + 10 (e^{0.2t}-1) \right),
\end{align*}
where $x_1^0 = 0.4 - \sqrt{0.2^2-0.1^2}\approx 0.227$ is the initial coordinate of the left contact point.

\begin{figure}[ht]
 \centering
 \subfigure[Boundary ELVIRA method.]{\includegraphics[width=0.45\columnwidth]{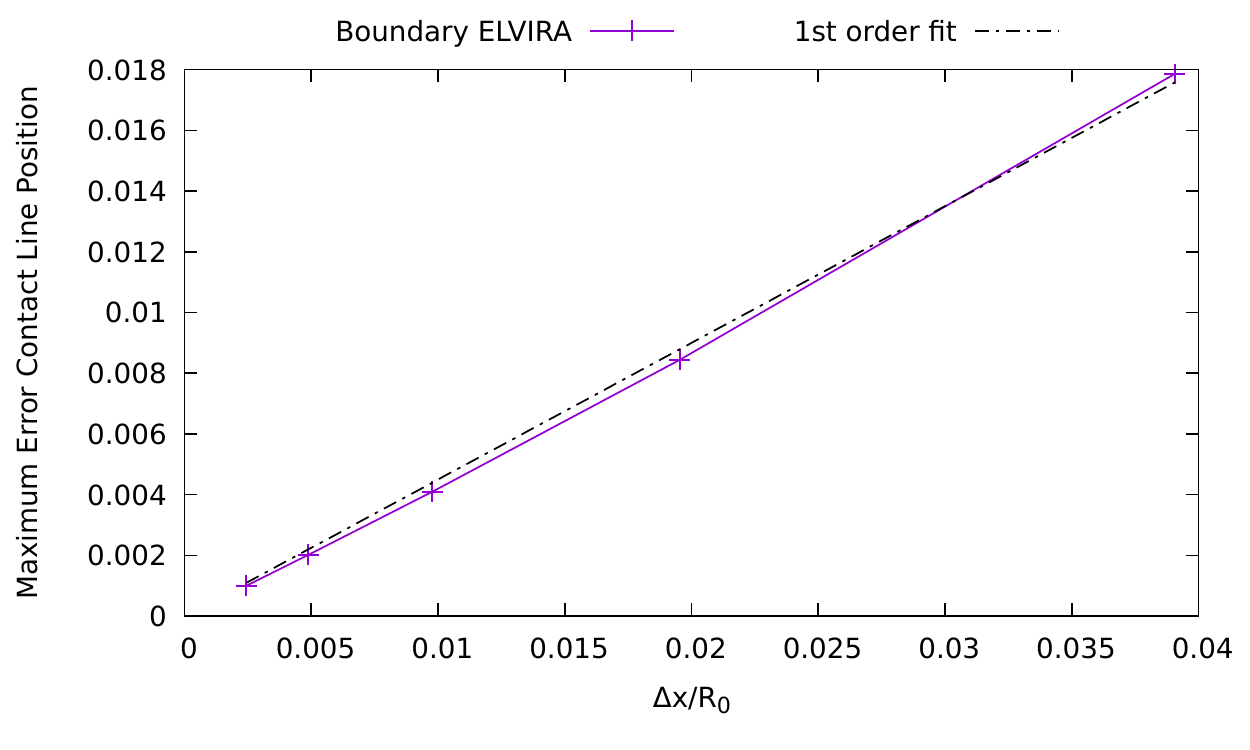}}
 \subfigure[Boundary Youngs method.]{\includegraphics[width=0.45\columnwidth]{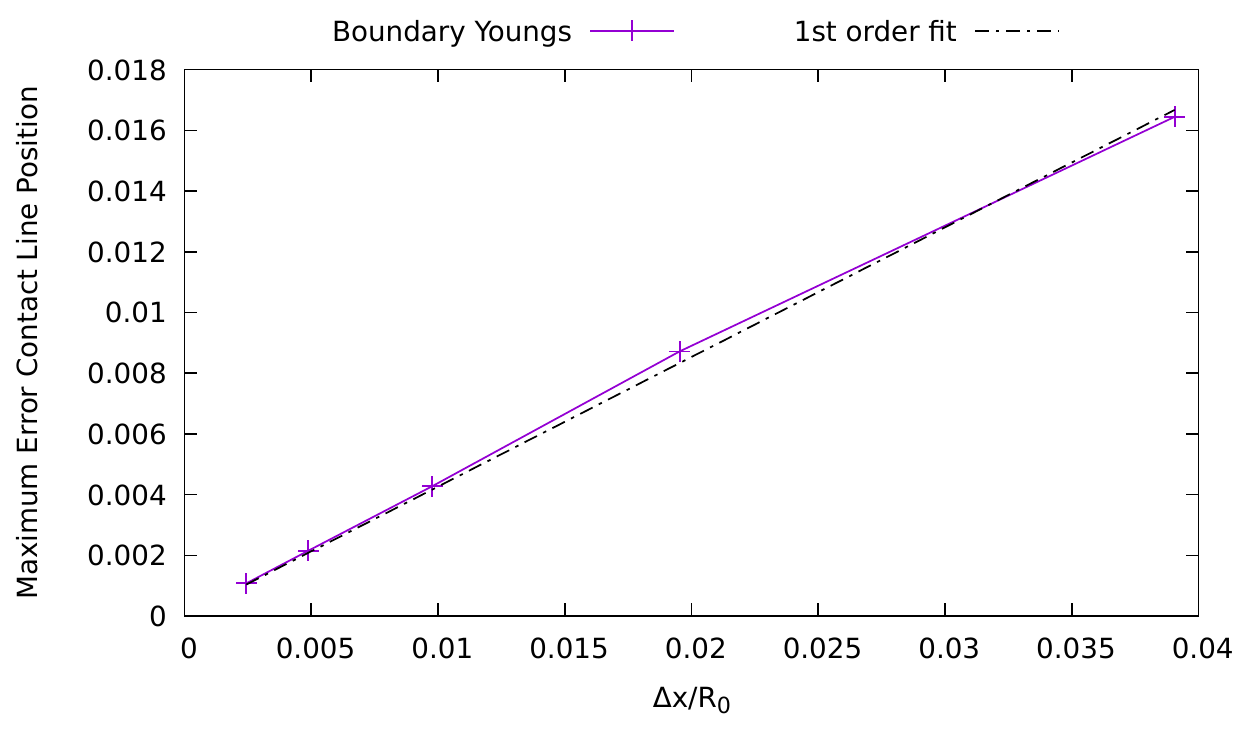}}
 \caption{Contact line evolution for the field \eqref{eqn:2D_navier_field}.}
 \label{fig:plot_navier_position_2d}
\end{figure}

\paragraph{Contact Line Motion:} Like in the previous example, both methods show first-order convergence with respect to the maximum norm regarding the motion of the contact line, see Figure~\ref{fig:plot_navier_position_2d}.

\begin{figure}[ht]
 \centering
 \subfigure[Boundary Youngs method.]{\includegraphics[width=0.45\columnwidth]{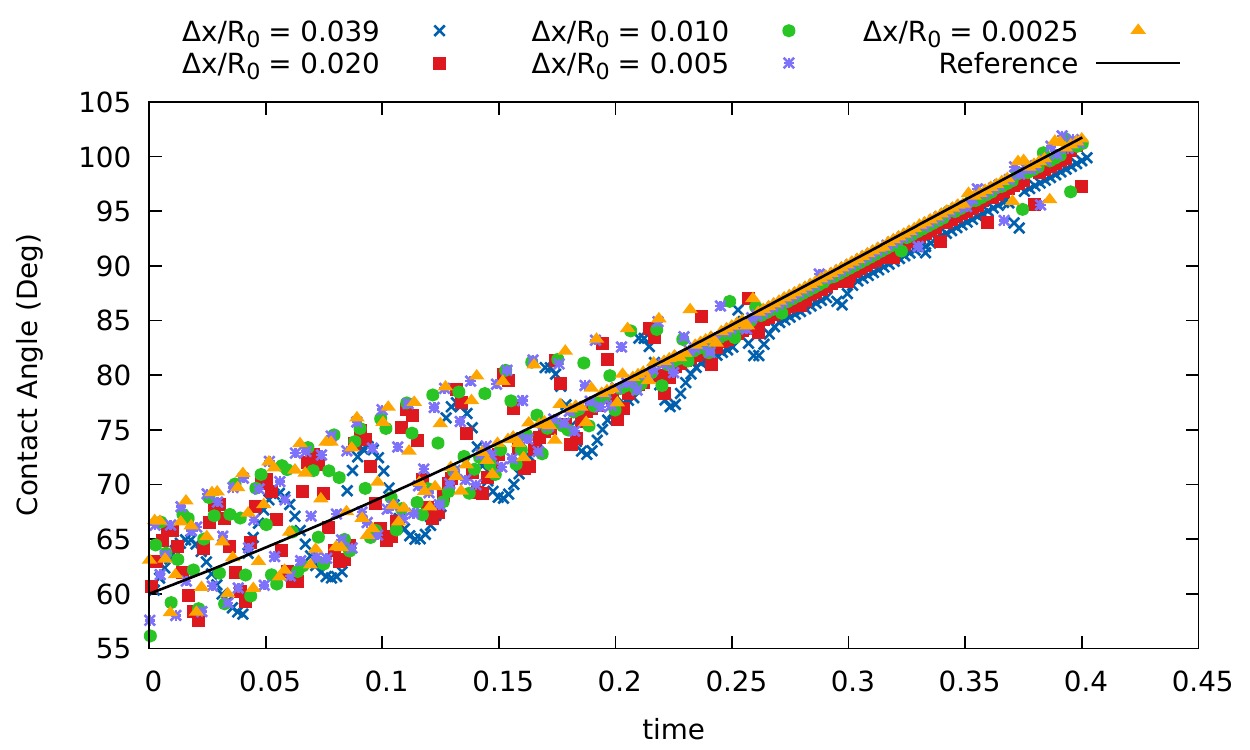}\label{fig:results_navier_youngs_2d}}
 \subfigure[Boundary ELVIRA method.]{\includegraphics[width=0.45\columnwidth]{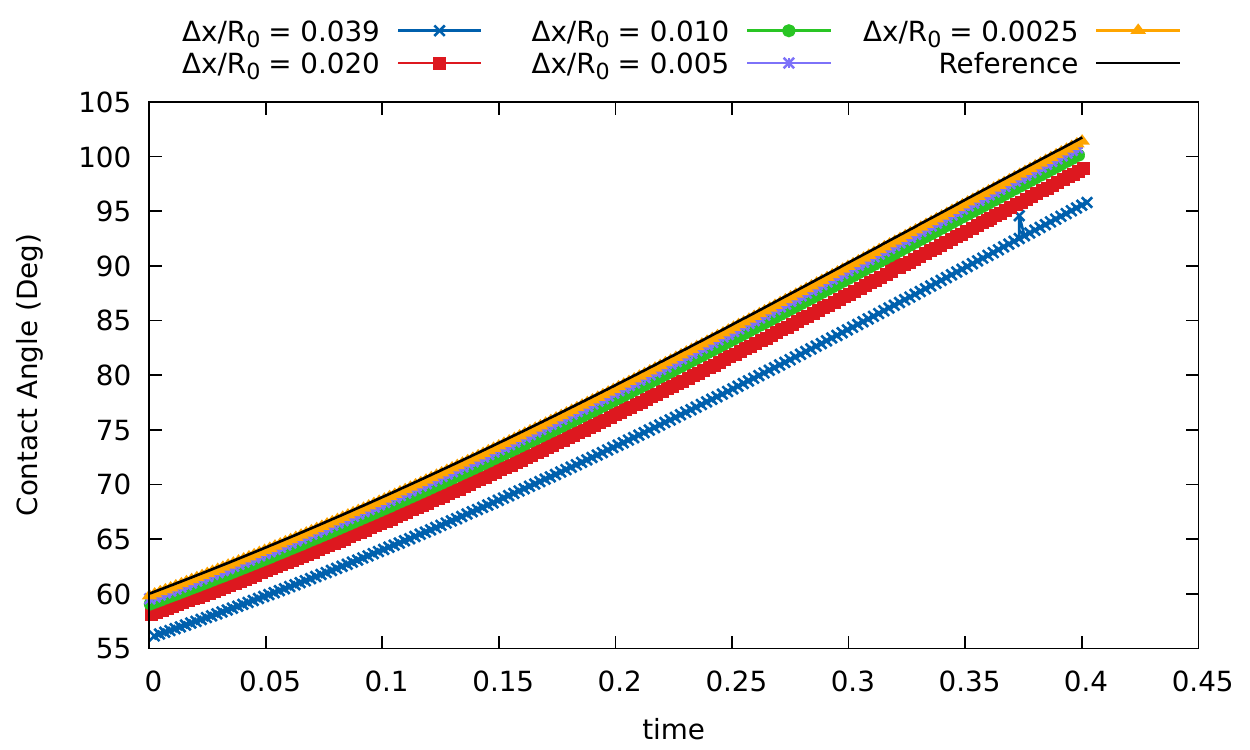}\label{fig:results_navier_elvira_2d}}
 \caption{Numerical contact angle evolution for the field \eqref{eqn:2D_navier_field}.}
 \label{fig:results_navier_2d}
\end{figure}

\paragraph{Contact Angle Evolution:} The numerical contact angle for the Boundary Youngs method is again subject to strong oscillations ($\pm 10^\circ$ in this case) and does not convergence with mesh refinement as visible in Figure~\ref{fig:results_navier_youngs_2d}. In contrast to that, the evolution of the contact angle for the Boundary ELVIRA method is first-order convergent and smooth with jumps visible only on a coarse grid, see Figure~\ref{fig:results_navier_elvira_2d}.

\begin{figure}[ht]
 \centering
 \subfigure[Boundary ELVIRA (for the contact angle).]{\includegraphics[width=0.45\columnwidth]{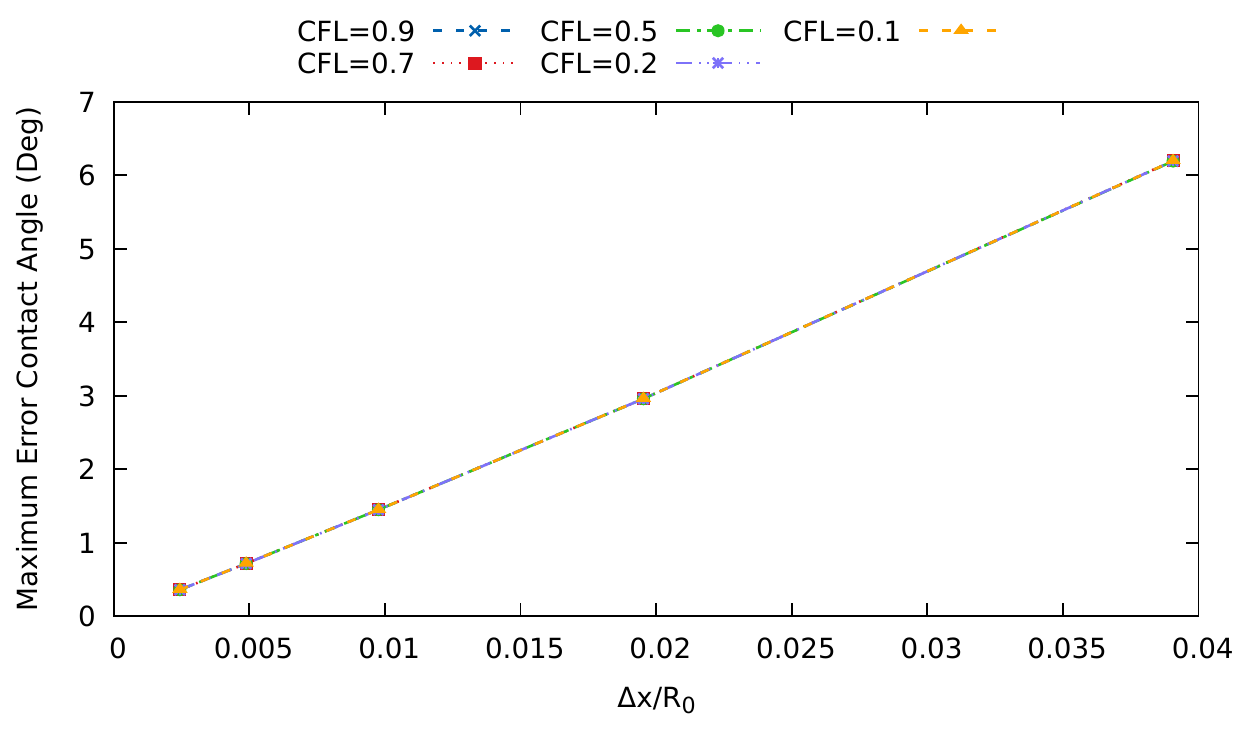}}
 \subfigure[Boundary Youngs (for the contact line position).]{\includegraphics[width=0.45\columnwidth]{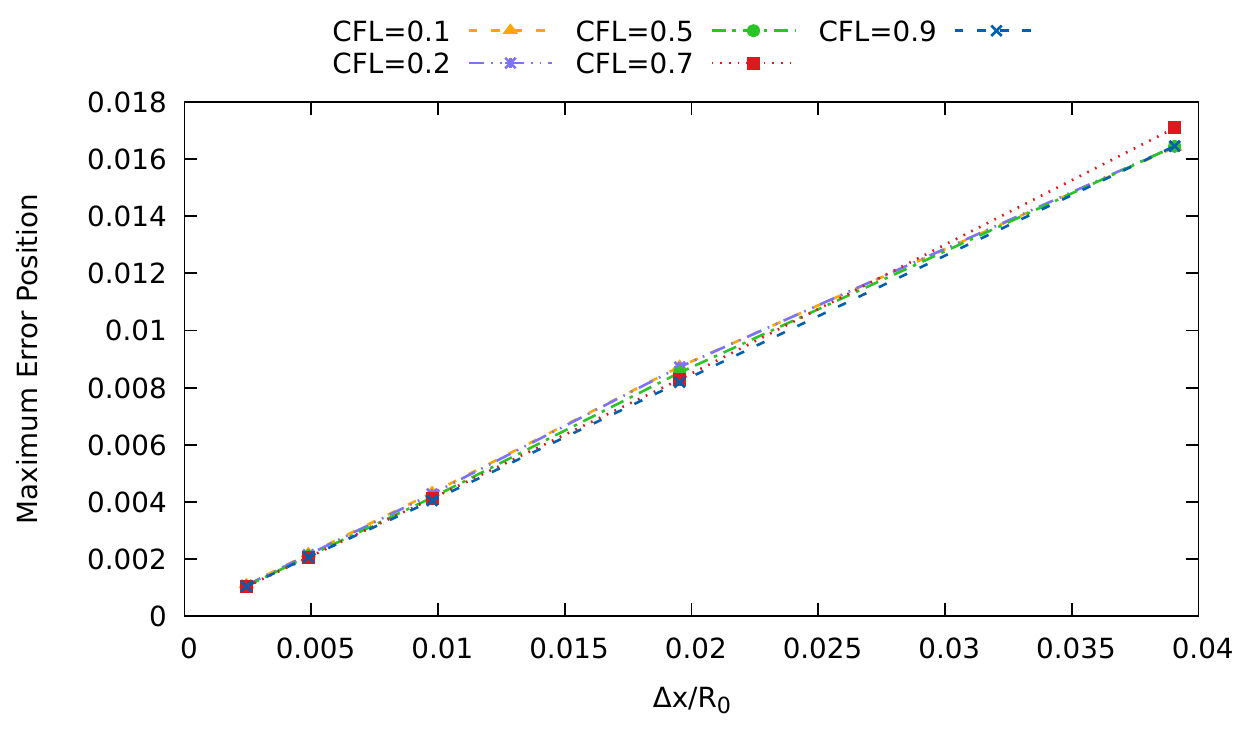}}
 \caption{Convergence behavior for the field \eqref{eqn:2D_navier_field} and different values of CFL.}
 \label{fig:cfl_study_1}
\end{figure}

\paragraph{Influence of the Courant number:} The influence of the Courant number for the convergence of the Boundary ELVIRA method (with respect to the contact angle) and the Boundary Youngs method (with respect to the contact line position) is shown in Figure~\ref{fig:cfl_study_1}. Apparently, there is hardly any influence in the considered example. Linear convergence is achieved for all reported Courant numbers from $0.1$ to $0.9$.\\
\\
Both numerical methods show excellent volume conservation. The relative volume error is at most of the order $10^{-10}$. Thanks to the volume redistribution algorithm, the volume fraction fields are exactly bounded up to machine precision. 

\subsection{Time-dependent Linear Velocity Field}
\label{section:time_dependent_fields}
As a third example, we consider the spatially linear, time-dependent velocity field of the form
\begin{align}
\label{eqn:time_dependent_field}
v(t,x_1,x_2) = \cos\left(\frac{\pi t}{\tau}\right) (v_0+c_1 x_1 + c_2 x_2, -c_1 x_2). 
\end{align}
As mentioned before, the time-dependent coefficient $\cos((\pi t)/\tau)$ is a classical choice to test advection methods by comparing the phase volumes at $t=0$ and $t=\tau$. Here we consider the full dynamics of the advection by solving the kinematic evolution equation for the field \eqref{eqn:time_dependent_field} explicitly.\\
\\
Using the ansatz \eqref{eqn:theta_ansatz}, it is easy to show that the exact solution for the latter velocity field is given by (for $c_1 \neq 0$)
\begin{equation}
\begin{aligned}
x_1(t) &= x_0 e^{c_1 s(t)} +\frac{v_0}{c_1}\left(e^{c_1 s(t)}-1 \right), \quad \thetaref(t) &= \frac{\pi}{2} + \arctan\left(-\cot \theta_0 e^{2 c_1 s(t)} \pm \frac{c_2}{2c_1}(e^{2 c_1 s(t)}-1) \right),
\end{aligned}
\end{equation}
where $s(t)$ is defined as
\[ s(t) = \frac{\tau \sin(\pi t/\tau)}{\pi}. \]
In particular, the evolution is periodic in $t$ with period $2\tau$. Note that the solution \eqref{eqn:exact_solution_2d} is recovered in the limit $\tau \rightarrow \infty$ since
\[ \lim_{\tau \rightarrow \infty} \frac{\tau \sin(\pi t/\tau)}{\pi} = \lim_{\tau \rightarrow \infty} \frac{\tau (\pi t/\tau)}{\pi} = t. \]
As a concrete example, we choose again $v_0 = -0.2$, $c_1 = 0.1$ and $c_2 = -2$ together with $\tau=0.2$.

\begin{figure}[ht]
 \centering
 \subfigure[Boundary Youngs Method.]{\includegraphics[width=0.45\columnwidth]{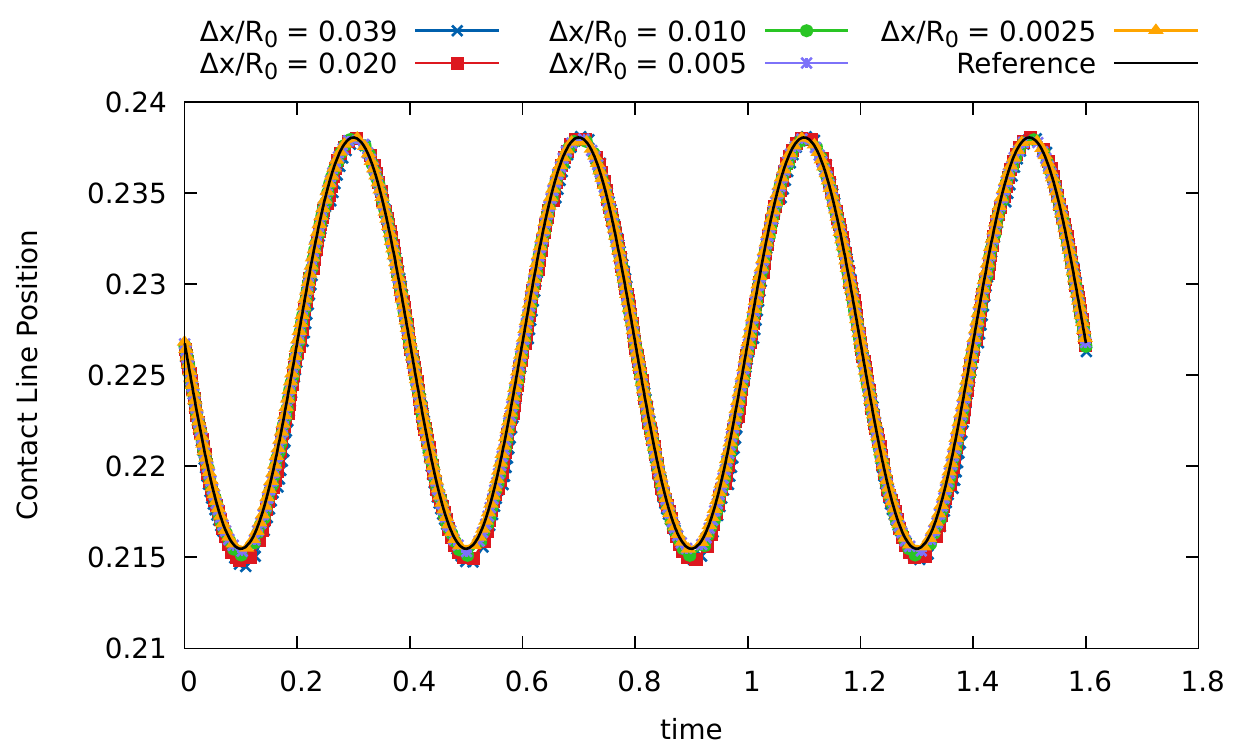}}
 \subfigure[Boundary ELVIRA Method.]{\includegraphics[width=0.45\columnwidth]{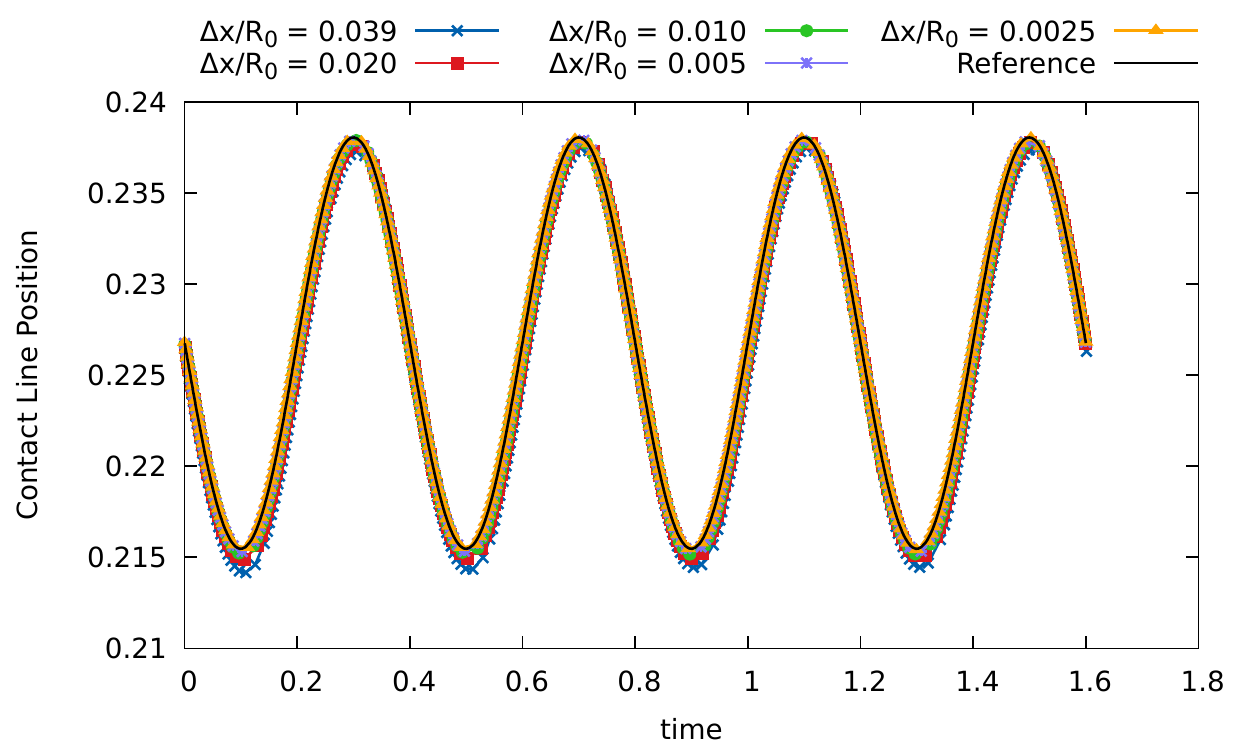}}
 \caption{Contact line motion for the field \eqref{eqn:time_dependent_field}.}
 \label{fig:plot_timedep_position}
\end{figure}

\paragraph{Contact Line Motion:} The contact line motion (see Figure~\ref{fig:plot_timedep_position}) is first-order convergent for both methods.

\begin{figure}[ht]
 \centering
 \subfigure[Boundary Youngs Method.]{\includegraphics[width=0.45\columnwidth]{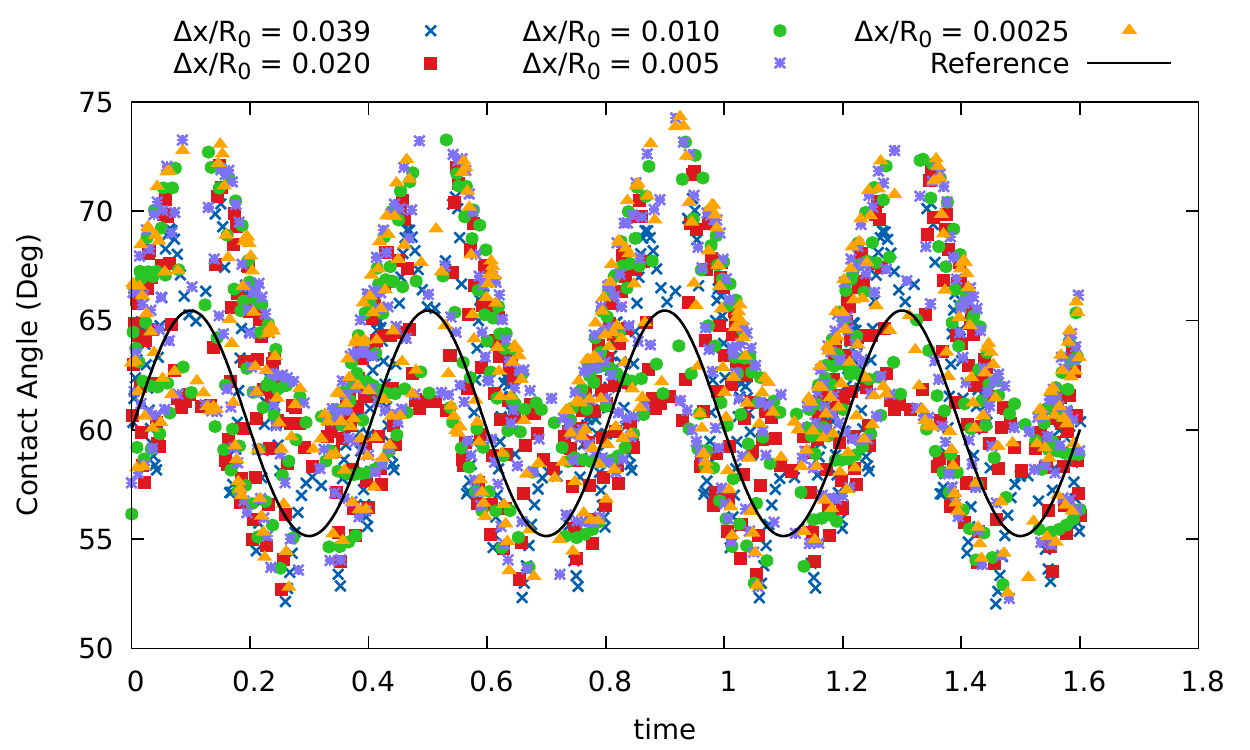}\label{fig:plot_timedep_contact_angle_youngs}}
 \subfigure[Boundary ELVIRA Method.]{\includegraphics[width=0.45\columnwidth]{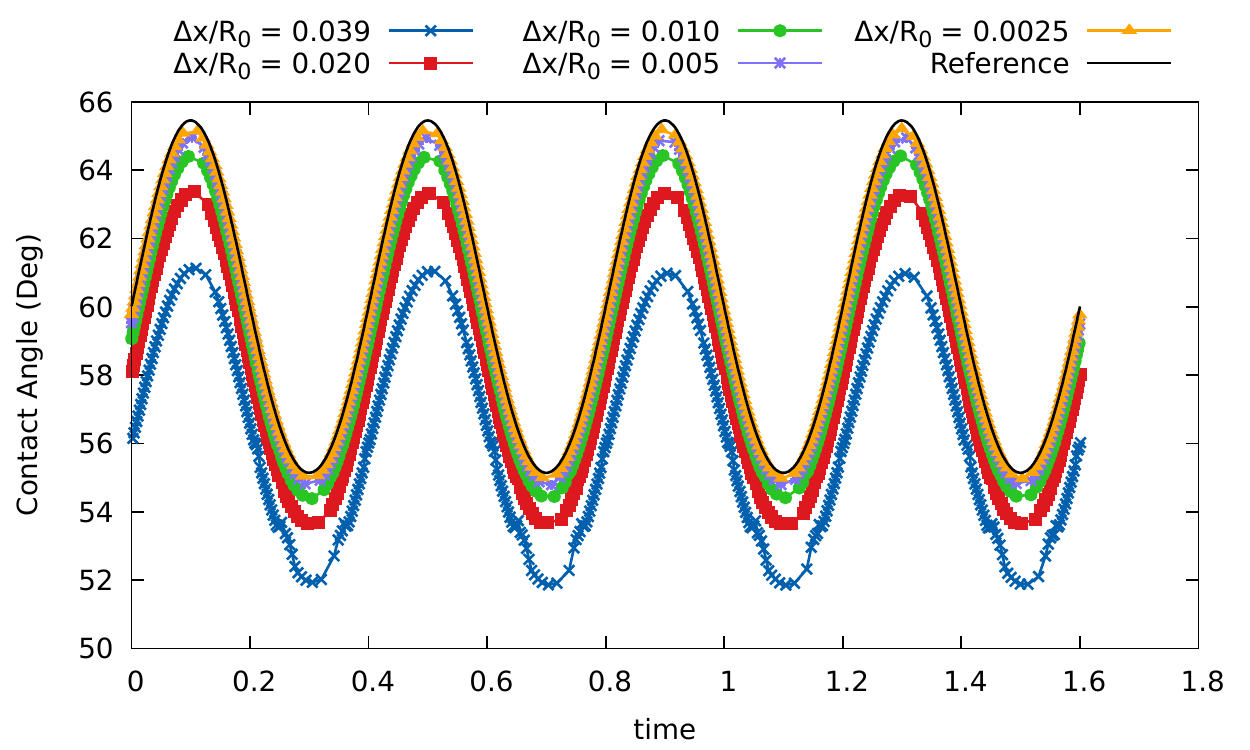}\label{fig:plot_timedep_contact_angle_elvira}}
 \caption{Contact angle evolution for the field \eqref{eqn:time_dependent_field}.}
 \label{fig:plot_timedep_contact_angle}
\end{figure}

\paragraph{Contact Angle Evolution:} Like in the previous examples, the numerical contact angle shows strong oscillations for the Boundary Youngs method, see Figure~\ref{fig:plot_timedep_contact_angle_youngs}. The evolution of the numerical contact angle for the Boundary ELVIRA method is reported in Figure~\ref{fig:plot_timedep_contact_angle_elvira}. While the numerical contact angle shows some deviations from the smooth reference curve on coarse grids, the period of the exact solution is still captured correctly. Like in the examples discussed before, refinement of the mesh at a fixed Courant number of $\cfl = 0.2$ leads to smoothening of the results and first-order convergence in the maximum norm. The results show the ability of the Boundary ELVIRA method to accurately capture the dynamics of the contact angle evolution.

\begin{figure}[ht]
 \centering
 \subfigure[Contact Angle.]{\includegraphics[width=0.45\columnwidth]{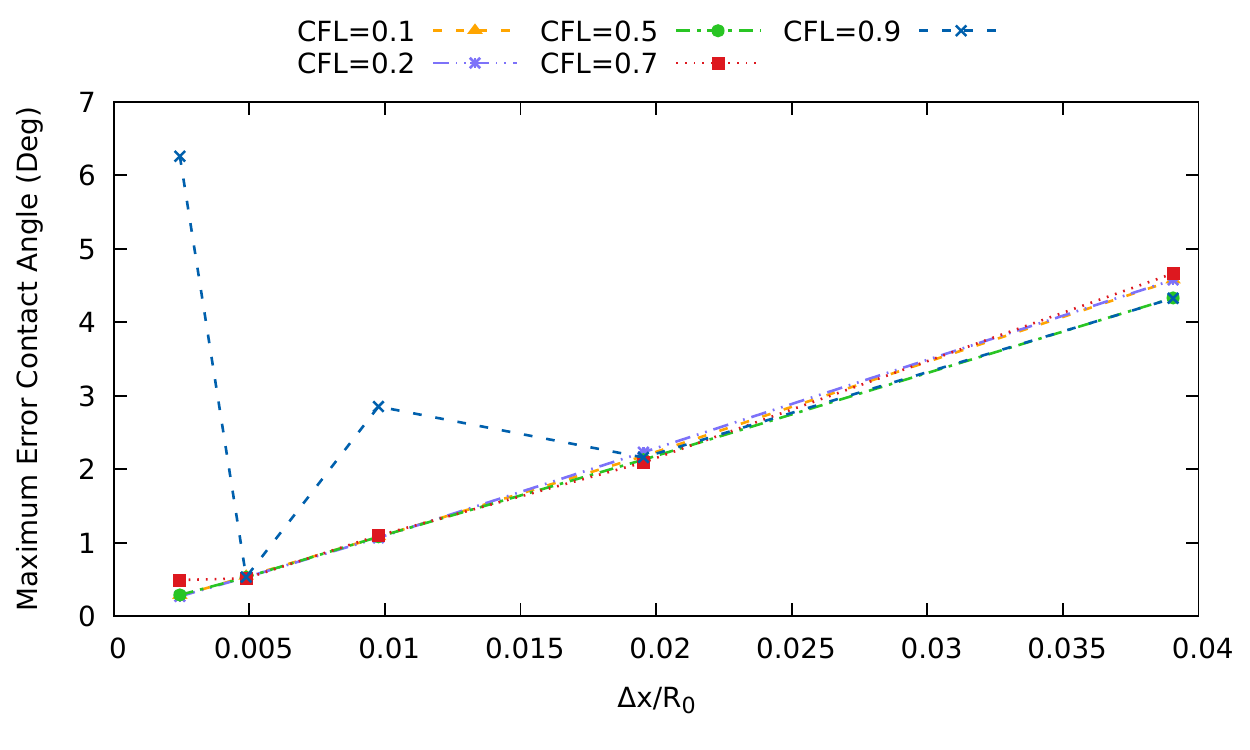}\label{fig:cfl_study_2_angle}}
 \subfigure[Contact Line Position.]{\includegraphics[width=0.45\columnwidth]{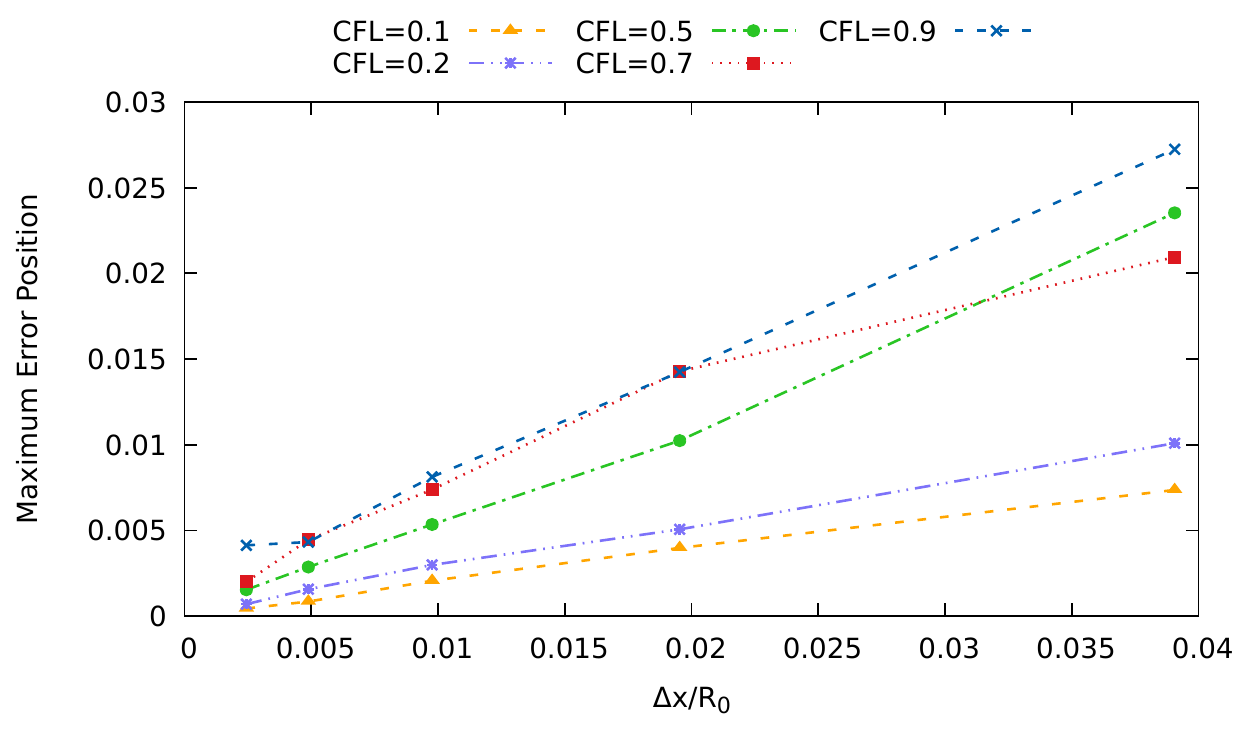}}
 \caption{Convergence behavior for the field \eqref{eqn:time_dependent_field} with the Boundary ELVIRA method.}
 \label{fig:cfl_study_2}
\end{figure}

\paragraph{Influence of the Courant number:} The results for the time-dependent linear field \eqref{eqn:time_dependent_field} turn out to be much more sensitive to the choice of the Courant number, see Figure~\ref{fig:cfl_study_2}.\\
The Boundary Youngs reconstruction fails completely for $\cfl$ greater than $0.5$ due to the appearance of interface cells with $|\nabla_h \alpha| \approx 0$. The convergence of the contact angle in the maximum norm for the Boundary ELVIRA method breaks down for $\cfl \geq 0.7$ (see Figure~\ref{fig:cfl_study_2_angle}).  Note, however, that a Courant number as large as $0.7$ is rarely achieved in multiphase flow simulations of systems governed by capillary effects (which is typically the case for wetting problems). In these systems, the numerical time step is usually limited by a stability criterion based on the propagation of capillary waves (see, e.g., \cite{Tryggvason2011}) and the (advective) CFL number is small.

\section{Summary and Conclusion}
The contact line advection problem is studied based on the geometrical Volume-of-Fluid method. Adaptations of the Youngs and ELVIRA methods to reconstruct the interface close to the boundary are introduced (see \cite{OnlineRepoVOF} for the implementations in FORTRAN). This allows to solve the transport equation for the interface without enforcing any boundary condition on the volume fraction field (except for inflow boundary conditions on artificial boundaries). Both the contact line position and the contact angle are evaluated based on the piecewise linear approximation of the interface (PLIC).\newline
\newline
The Boundary Youngs method allows to track the motion of the contact line with first-order accuracy. However, a meaningful evaluation of the contact angle in terms of the local interface orientation is not possible. Instead, the numerical contact angle shows strong oscillations resulting from a spatial dependence of the reconstruction error which is already present for a planar interface.\newline
\newline
The Boundary ELVIRA method delivers first-order convergent results for the dynamics of both the contact line motion and the contact angle evolution. The results are verified using an explicit and a numerical solution of the kinematic evolution equation \eqref{eqn:theta_evolution_equation}.\newline
\newline
Based on the Boundary ELVIRA method to transport the contact angle, one may develop numerical methods of dynamic wetting, where the contact angle is not prescribed as a geometric boundary condition. Instead, a local force term depending on the interface orientation may be introduced. The latter approach, as discussed e.g. in \cite{Sikalo2005} for an algebraic VOF method, avoids the necessity to manually ``adjust'' the contact angle after a transport step. The present method allows to evaluate the local interface orientation with high accuracy in geometrical Volume of Fluid methods and may, therefore, serve as a basis for future numerical methods in dynamic wetting.

\paragraph{Acknowledgments:} We kindly acknowledge the financial support by the German Research Foundation (DFG) within the Collaborative Research Centre 1194 “Interaction of Transport and Wetting Processes”, Projects B01 and Z-INF. Calculations for this research were conducted on the Lichtenberg high performance computer of the TU Darmstadt.

\bibliography{ms.bbl}

\clearpage
\appendix
\section[Appendix]{Appendix}
The following definition of a $\mathcal{C}^{1,2}$-family of moving hypersurfaces can also be found in \cite{Kimura.2008}, \cite{Pruss.2016} and in a similar form in \cite{Giga.2006}.
\begin{definition}
Let $I=(a,b)$ be an open interval. A family $\{\Sigma(t)\}_{t \in I}$ with $\Sigma(t) \subset \RR^3$ is called a \emph{$\mathcal{C}^{1,2}$-family of moving hypersurfaces} if the following holds.
\begin{enumerate}[(i)]
 \item Each $\Sigma(t)$ is an orientable $\mathcal{C}^2$-hypersurface in $\RR^3$ with unit normal field denoted as $\nsigma(t,\cdot)$.
 \item The graph of $\Sigma$, given as
 \begin{align}
 \label{eqn:def_moving_interface}
 \move := \gr \Sigma = \bigcup_{t \in I} \{t\} \times \Sigma(t) \subset \RR\times\RR^{3},
 \end{align}
 is a $\mathcal{C}^1$-hypersurface in $\RR \times \RR^3$.
 \item The unit normal field is continuously differentiable on $\move$, i.e.
 \[ \nsigma \in \mathcal{C}^1(\move). \]
\end{enumerate}
A family $\{\overline{\Sigma}(t)\}_{t \in I}$ is called a \emph{$\mathcal{C}^{1,2}$-family of moving hypersurfaces with boundary $\partial\Sigma(t)$} if the following holds.
\begin{enumerate}[(i)]
 \item Each $\overline{\Sigma}(t)$ is an orientable $\mathcal{C}^2$-hypersurface in $\RR^3$ with interior $\Sigma(t)$ and non-empty boundary $\partial\Sigma(t)$, where the unit normal field is denoted by $\nsigma(t,\cdot)$.
 \item The graph of $\overline{\Sigma}$, i.e.
 \[ \gr \overline{\Sigma} = \bigcup_{t \in I} \{t\} \times \overline{\Sigma}(t) \subset \RR\times\RR^{3}, \]
 is a $\mathcal{C}^1$-hypersurface with boundary $\gr(\partial\Sigma)$ in $\RR\times\RR^3$.
  \item The unit normal field is continuously differentiable on $\gr \overline{\Sigma}$, i.e.
 \[ \nsigma \in \mathcal{C}^1(\gr \overline{\Sigma}). \]
\end{enumerate}
\end{definition}
Being the boundary of a submanifold with boundary, the set $\gr(\partial\Sigma)$ is itself a submanifold (without boundary).

\begin{remark}[Construction of an analytical solution]
Equation \eqref{eqn:theta_evolution_linear_2d} may be solved with the following Ansatz: We look for solutions of the form
\begin{align}
\label{eqn:theta_ansatz}
\theta(t) = \frac{\pi}{2}+\arctan(f(t)). 
\end{align}
It is an easy exercise to show that this yields the following ordinary differential equation for $f$:
\begin{align}
\label{eqn:transformed_evolution_problem}
\dot{f} = \pm c_2 + 2 c_1 f. 
\end{align}
The initial condition for $\theta$ translates to
\begin{align}
\label{eqn:transformed_evolution_problem_initial_condition}
f(0) = -\cot \theta_0. 
\end{align}
For $c_1 \neq 0$, the initial-value problem \eqref{eqn:transformed_evolution_problem}-\eqref{eqn:transformed_evolution_problem_initial_condition} has the unique solution
\[ f(t) = -\cot \theta_0 \, e^{2c_1 t} \pm c_2 \frac{e^{2c_1 t} -1}{2c_1}. \]
Hence, we obtain the desired solution \eqref{eqn:exact_solution_2d}.
\end{remark}


\end{document}